\RequirePackage{fix-cm}
\makeatletter
\def\cl@chapter{}
\makeatother
\documentclass[smallextended]{svjour3}       

\smartqed  
\usepackage{graphicx}
\usepackage{mathptmx}      
\usepackage{bmkstyle}
\usepackage{bmkstyle_maths}
\usepackage[capitalise,nameinlink]{cleveref}
\crefname{section}{Sect.}{Sect.}
\Crefname{section}{Section}{Sections}
\crefname{listing}{\lstlistingname}{\lstlistingname}
\Crefname{listing}{Listing}{Listings}

\newtheorem{assumption}{Assumption}

\definecolor{otherblue}{rgb}{0,0,0}
\def\rbl#1{#1}

\journalname{Springer }

\begin{document}
	\title{\rbl{Efficient} Adaptive Stochastic Collocation Strategies for Advection-Diffusion Problems with Uncertain Inputs
		\thanks{{This work was partially supported by the EPSRC grants EP/V048376/1 and EP/W033801/1.}}
	}
	
	\titlerunning{\rbl{Efficient} Adaptive SC-FEM for Advection--Diffusion}
	
	\author{Benjamin M. Kent         \and
		Catherine E. Powell	\and
		David J. Silvester		\and
		Małgorzata J. Zimoń 
	}
	
	
	\institute{B.M. Kent \at
		Department of Mathematics, University of Manchester, Manchester, M13 9PL, UK
		\email{benjamin.kent@mancehster.ac.uk}
		\and
		C.E. Powell \at
		Department of Mathematics, University of Manchester, Manchester, M13 9PL, UK
		\email{catherine.powell@manchester.ac.uk}
		\and
		D.J. Silvester \at
		Department of Mathematics, University of Manchester, Manchester, M13 9PL, UK
		\email{david.silvester@manchester.ac.uk}
		\and
		M.J. Zimoń \at
		IBM Research Europe,
		Daresbury Laboratory,
		Keckwick Lane, Warrington,
		WA4 4AD,
		UK
		\email{malgorzata.zimon@uk.ibm.com}
	}
	
	\date{Received: date / Accepted: date}

	\maketitle
	
	\begin{abstract}
		Physical models with uncertain inputs are commonly represented as parametric partial differential equations (PDEs). That is, PDEs with inputs that are expressed as functions of parameters with an associated probability distribution. Developing efficient and accurate solution strategies that account for errors on the space, time and parameter domains simultaneously is highly challenging. Indeed, it is well known that standard polynomial-based approximations on the parameter domain can incur errors that grow in time. In this work, we focus on advection–diffusion problems with parameter-dependent wind fields. A novel adaptive solution strategy is proposed that allows users to combine stochastic collocation on the parameter domain with off-the-shelf adaptive timestepping algorithms with local error control. This is a non-intrusive strategy that builds a polynomial-based surrogate that is adapted sequentially in time. The algorithm is driven by a so-called hierarchical estimator for the parametric error and balances this against an estimate for the global timestepping error which is derived from a scaling argument. 
			
		\keywords{parametric PDEs, stochastic collocation, adaptivity, error estimation}
		\subclass{65C20 \and 65M22 \and 65M15 \and 65M70}
	\end{abstract}
	
	\section{Introduction}
	\label{sec:intro}
	Parametric partial differential equations (PDEs) are commonly used  to model physical problems with uncertain inputs that are represented as random variables. Briefly, if the uncertain inputs can be expressed as functions of a finite set of independent random variables $\xi_{i}: \Omega \to \Gamma_{i} \subset \R$, $i=1,\ldots, d$ on a probability space $(\Omega, {\cal F}, \mathbb{P})$, then with the change of variables $y_{i}=\xi_{i}(\omega)$, $i=1,\ldots d$, we may recast the inputs and PDE solution as functions of the parameter vector $\y= [y_{1}, \ldots, y_{d}]^{\top} \in \Gamma$ where $\Gamma:= \Gamma_{1} \times \cdots \times \Gamma_{\nY} \subset \R^\nY$ denotes the parameter domain.  Denoting the Borel $\sigma$-algebra by $\mathcal{B}(\Gamma)$ and the joint probability density by $\rho(\y) = \prod_{\iY=1}^{\nY} \rho_{\iY} (y_{\iY})$, where  $\rho_{\iY} (y_{\iY})$ is the probability density of $\xi_{i}$ we can formulate the problem and compute statistical quantities of interest with respect to the equivalent probability space $(\Gamma, \mathcal{B}(\Gamma), \rho(\y) d\y)$. 
	
	In this work, we are interested in the efficient and accurate numerical solution of parametric time-dependent advection-diffusion problems of the form
	\begin{equation}\label{first_equation}
	\ddt u(\x,t,\y) - \diffcnst \Delta u(\x,t,\y)  + \wind(\x,\y) \cdot \nabla u(\x,t,\y) = 0,
	\end{equation}
	where the $d$ input parameters $y_{i}$ are images of random variables that encode our uncertainty about the spatially-varying advection field or `wind' $\wind$. Such time-dependent parametric PDEs are more challenging to solve than their steady-state counterparts. Different choices of $\y$ yield different space-time solutions that may be pushed further away from each other in time. This makes $u(\x,t,\y)$ harder to approximate as a function of $\y$ as time $t\to \infty$. This is not a new observation, for example it is discussed in \cite{Wan2006}, where the authors investigate a transport problem with uncertain velocity.

	Numerical solution strategies for parametric PDEs can be categorised as either non-intrusive methods, which involve the solution of decoupled deterministic problems for particular choices of $\y$ or intrusive methods, based on Galerkin projection, which require the solution of a single large coupled problem.  Among these, polynomial-based approximation schemes such as stochastic collocation  \cite{Xiu2005,Nobile2008a,Babuska2007} and generalised polynomial chaos (gPC) approximation \cite{Xiu2002,Babuska2004} that lead to surrogates that can be cheaply evaluated for any $\y \in \Gamma$ are useful for tasks like uncertainty quantification, optimisation and design. Non-intrusive schemes that allow the reuse of existing space-time solvers are especially appealing for practitioners. 
	
	
	
	Both intrusive and non-intrusive solution strategies have been extensively studied for parametric elliptic PDEs and state-of-the-art adaptive solvers exist for such problems.
	Anisotropic \cite{Nobile2008b} and adaptive \cite{Nobile2016} collocation schemes have been proposed for problems with a moderately large number of parameters.
	Adaptive variants of sparse grid collocation typically use a strategy based on Gerstner and Griebel's quadrature algorithm \cite{Gerstner2003}, creating an operator that is a sum of appropriately chosen hierarchical surplus operators.
	More recently, residual based a posteriori error estimation strategies have been developed for stochastic collocation approximations of solutions to parametric elliptic PDEs \cite{Guignard2018} and the convergence of adaptive algorithms driven by such estimators has been analysed \cite{Eigel2022,Feischl2021}.
	Less work has been done for time-dependent parametric PDEs.
	However, a residual-based error estimator was used to drive an adaptive stochastic collocation algorithm for the heat equation in \cite{Nobile2021} and time-adapted bases were used in the context of dynamical low rank approximation for parametric parabolic problems in \cite{Kazashi2021}.
	For gPC approximations, the use of polynomial bases that are adapted in time is discussed in \cite{Gerritsma2010}. 
	
	We will use sparse grid stochastic collocation \cite{Smolyak1963,Barthelmann2000,Novak1999} for the parametric approximation of \eqref{first_equation}. 
	This is appealing due to its non-intrusive nature and the fact that it leads to a polynomial surrogate.
	The main challenge we address is how to combine such a method with an adaptive timestepping scheme so that the distinct contributions to the total error remain balanced.
	For this, an efficient a posteriori error estimator is required.
	{\color{otherblue} The strategy we propose is different to the one proposed in \cite{Nobile2021}.
	Our estimator is not residual-based.
	Moreover, we use adaptive timestepping methods in a black-box fashion and adapt the polynomial approximation \emph{sequentially} in time.}
	
	
	
	%

	
	In \cref{sec:problem-statement} we introduce the parametric problem of interest in more detail and describe the numerical solution approach, combining sparse grid stochastic collocation on $\Gamma$ and an adaptive timestepping algorithm with local error control.
	In \cref{sec:motivational-example}, we consider a related ODE test problem and motivate the need to adapt the parametric approximation in time.
	An a posteriori error estimator combining a Gerstner and Griebel type strategy on the parameter domain and a global timestepping error estimation strategy is defined in \cref{sec:error-est}.
	A novel adaptive algorithm that is driven by the error estimator is described in \cref{sec:error-control} and numerical results are presented in \cref{sec:numerical-results} {and \cref{sub:example-sixtyfour}}. 
	
	\section{{Problem Statement and Approximation Scheme}}
	\label{sec:problem-statement}
	
	
	Let $\spatialdomain \subset \R^{n}$ ($n=2,3$) be an open polygonal spatial domain with Lipschitz boundary $\spatialbdry$ and let $\y= [y_{1}, \ldots, y_{d}]^{\top}$ be a vector of $d$ parameters taking values in a parameter domain $\Gamma= \Gamma_{1} \times \cdots \times \Gamma_{\nY} \subset \R^\nY$, with associated joint probability density $\rho(\y) = \prod_{\iY=1}^{\nY} \rho_{\iY} (y_{\iY})$.
	We consider an advection-diffusion problem with a fixed diffusion coefficient $\diffcnst > 0$ and a parameter-dependent, spatially divergence-free advection field $\wind:\spatialdomain \times \Gamma \to \R^\nX$.
	That is, we seek a solution $u:\spatialdomain \times [0,T] \times \Gamma \to \R$ satisfying
	\begin{equation}
	\begin{cases}
	\begin{aligned}
	\ddt u(\x,t,\y) &- \diffcnst \Delta u(\x,t,\y) + \wind(\x,\y) \cdot \nabla u(\x,t,\y) = 0, \end{aligned} & (\x,t) \in \spatialdomain \times (0,T],\\
	u(\x,0,\y) = u_0(\x), & \x \in \spatialdomain,
	\end{cases}
	\label{eq:adv-diff}
	\end{equation}
	$\rho$-almost surely on $\Gamma$ subject to the Dirichlet boundary condition
	\begin{equation}
	u(\x,t,\y) = u_{\spatialbdry}(\x,t), \quad (\x,t,\y) \in \spatialbdry \times [0,T]\times \Gamma.
	\end{equation}
	
	For the weak formulation, we define $\HOneBdry:=\{v\in \HOne \st \trace v = u_{\spatialbdry}\}$ and seek $u(\cdot,t,\y) \in \HOneBdry$ with $\frac{\partial u}{\partial t}(\cdot,t,\y) \in \HMinus$ that solves \begin{equation}
	\begin{cases}
	\begin{aligned}\int_{\spatialdomain} \ddt u(&\x,t,\y) v(\x)+ \diffcnst \nabla u(\x,t,\y) {\color{red} \cdot} \nabla v(\x) \\ &+ \wind(\x,\y) \cdot \nabla u(\x,t,\y) v(\x) \dd \x= 0, \end{aligned} & t \in (0,T],\\ \\
	\int_{\spatialdomain} u(\x,0,\y) v(\x) \dd \x= \int_{\spatialdomain} u_0(\x) v(\x) \dd \x,
	\end{cases}
	\label{eq:adv-diff-weak}
	\end{equation}
	$\rho$-almost surely on $\Gamma$ for all test functions $v \in \HOneZero:=\{v\in\HOne \st \trace v = 0 \}$.
	Conditions for existence and uniqueness of a solution for a fixed $\y$ are outlined in Evans \cite[Section 7.1.2]{Evans2010}.
	In particular, we need to assume that $\diffcnst>0$ is finite and $\wind(\x, \y)=[w^1(\x,\y), w^2(\x,\y)]^{\top}$ where $w^1,w^2 \in L^\infty(\spatialdomain)$, $\rho$-almost surely on $\Gamma$.
	In our numerical experiments we will choose wind fields of the form
	\begin{equation}
	\wind(\x,\y) := \wind_0(\x) + \sum_{\iY=1}^{\nY} y_{\iY} \wind_{\iY}(\x)
	\label{eq:wind-lincomb}
	\end{equation}
	where $\{\wind_{\iY}\}_{\iY=0}^{\nY}$ are spatial fields with $w^1_{\iY},w^2_{\iY} \in L^\infty(\spatialdomain)$ for each $\iY=0,...,\nY$ and the parameters $y_{i} \in [-1,1]$ (i.e., the associated random variables are {bounded}).
	
	
	

	For the spatial discretisation, we apply finite element approximation. We do not specify a particular method here but assume that the chosen approximation is globally continuous and defined piecewise with respect to a spatial mesh of non-overlapping triangular (tetrahedral if $n=3$) or rectangular (hexahedral) elements with characteristic edge length $\h$. The associated nodes are partitioned into two sets $\nodes=\{\x_\iSpatial\}_{\iSpatial=1}^{\nInterior} \cup \{\x_\iSpatial\}_{\iSpatial=\nInterior+1}^{\nInterior+\nBoundary}$ corresponding to interior and boundary nodes, respectively, as are the associated sets of nodal basis functions $\{ \spatialbasis_{\iSpatial} \}_{\iSpatial=1}^{\nInterior} \cup \{ \spatialbasis_{\iSpatial} \}_{\iSpatial=\nInterior + 1}^{\nInterior+\nBoundary}$. For each $(t,\y) \in [0,T] \times \Gamma$ we seek an approximation of the form 
	\begin{equation}
		u^{\h}(\x,t,\y)=\sum_{\iSpatial=1}^{\nInterior} u^{\h}_\iSpatial(t,\y) \spatialbasis_\iSpatial(\x) + \underbrace{\sum_{\iSpatial=1}^{\nBoundary} \udirichlet(\x_{\nInterior+\iSpatial},t) \spatialbasis_{\nInterior+\iSpatial} (\x)}_{=:\udirichlet^{\h}(\x,t)}
	\end{equation}
	where the coefficients associated with the boundary nodes are fixed. Substituting this into \eqref{eq:adv-diff-weak} and using the test space $\testSpace=\{\spatialbasis_\iSpatial\}_{\iSpatial=1}^{\nInterior} \subset \HOneZero$ yields a parametric system of ODEs for the vector of coefficients $\vec{u}^\h (t,\y):=[u^{\h}_1(t,\y),...,u^{\h}_{\nInterior}(t,\y)]^{\top}$. To formulate this system, we define a mass matrix $\massMatrix$, a diffusion matrix $\diffMatrix$ and a parameter-dependent advection matrix $\advMatrix(\y)$ associated with the interior nodes by
	\begin{equation}
	\begin{aligned}
	\left[\massMatrix\right]_{\iSpatial,\iSpatialAlt} &= \int_\spatialdomain \spatialbasis_{\iSpatial}(\x) \spatialbasis_{\iSpatialAlt}(\x) \dd \x,\\ \left[\diffMatrix\right]_{\iSpatial,\iSpatialAlt} &= \int_\spatialdomain \nabla\spatialbasis_{\iSpatial}(\x) \cdot \nabla\spatialbasis_{\iSpatialAlt}(\x) \dd \x,\\ \left[\advMatrix(\y)\right]_{\iSpatial,\iSpatialAlt} &= \int_\spatialdomain \left(\wind(\x,\y) \cdot \nabla \spatialbasis_{\iSpatialAlt}(\x)\right) {\spatialbasis_{\iSpatial}(\x)} \dd \x
	\end{aligned} \label{eq:fem-matrices}
	\end{equation}
	for $\iSpatial,\iSpatialAlt = 1,...,\nInterior$. We also define analogous matrices $\massMatrix_{\spatialbdry}, \diffMatrix_{\spatialbdry}$ and $\advMatrix_{\spatialbdry}(\y)$  that account for boundary interactions, for $\iSpatial = 1,..., \nInterior$ and $\iSpatialAlt=1,...,\nBoundary$, as follows
	
	\begin{equation}
	\begin{aligned}
	\left[\massMatrix_{\spatialbdry}\right]_{\iSpatial,\iSpatialAlt} &= \int_\spatialdomain \spatialbasis_{\iSpatial}(\x) \spatialbasis_{\nInterior+\iSpatialAlt}(\x) \dd \x,\\ \left[\diffMatrix_{\spatialbdry}\right]_{\iSpatial,\iSpatialAlt} &= \int_\spatialdomain \nabla\spatialbasis_{\iSpatial}(\x) \cdot \nabla\spatialbasis_{\nInterior+\iSpatialAlt}(\x) \dd \x,\\ \left[\advMatrix_{\spatialbdry}(\y)\right]_{\iSpatial,\iSpatialAlt} &= \int_\spatialdomain \left(\wind(\x,\y) \cdot \nabla \spatialbasis_{\nInterior+\iSpatialAlt}(\x)\right) { \spatialbasis_{\iSpatial}(\x) } \dd \x.
	\end{aligned}
	\end{equation} 
	
	The coefficient vector $\vec{u}^\h:[0,T]\times \Gamma \to \R^{\nInterior}$ then satisfies the ODE system
	\begin{equation}
	\begin{cases}
	\massMatrix \ddt \vec{u}^\h(t,\y) + \diffcnst \diffMatrix \vec{u}^\h(t,\y) + \advMatrix(\y) \vec{u}^\h(t,\y) = \forceFn(t,\y), & t \in (0,T],\\
	\vec{u}^\h(0,\y) = \vec{u}^\h_0=[u_{0}(\x_{1}), ..., u_{0}(\x_{\nInterior})]^\top,&
	\end{cases}
	\label{eq:system-of-odes}
	\end{equation}
	$\rho$-almost surely on $\Gamma$ subject to the forcing term
	\begin{equation}
	\forceFn_{}(t,\y):= - \massMatrix_{\spatialbdry} \frac{ \dd \vec{u}^\h_{\spatialbdry} (t)}{\dd t} - \left(\diffcnst  \diffMatrix_{\spatialbdry}  + \advMatrix_{\spatialbdry}(\y) \right) \vec{u}^\h_{\spatialbdry}(t), \label{eq:fem-f}
	\end{equation}
	where $\vec{u}^\h_{\spatialbdry}(t) := [u_{\spatialbdry}(\x_{\nInterior +1},t), ..., u_{\spatialbdry}(\x_{\nInterior +\nBoundary},t)]^\top$. After solving this system, the spatially-discrete approximation $u^{\h}:\spatialdomain \times [0,T] \times \Gamma \to \R$ can be reconstructed as
	\begin{equation}
	\begin{aligned}
	u^\h(\x,t,\y) := \sum_{\iSpatial=1}^{\nInterior} \left[\vec{u}^\h(t,\y)\right]_\iSpatial \spatialbasis_\iSpatial(\x) + \udirichlet^{\h}(\x,t).
	\end{aligned}
	\label{eq:full-functional-approximation}
	\end{equation}
	Note that in general, $\udirichlet^{\h} \neq \udirichlet$ unless the boundary condition is a piecewise polynomial of an appropriate degree. However, we will ignore this error here. Our focus is on approximating the solution of the ODE system \eqref{eq:system-of-odes}, or equivalently the function $u^{h}$ in \eqref{eq:full-functional-approximation}, that arises after applying the spatial discretisation.
	
	
	
	

	
	
	For the parametric approximation, we will construct a global polynomial approximation on $\Gamma$ that interpolates $u^{\h}$ at a set of collocation points. Specifically, we use the sparse grid construction proposed by Smolyak \cite{Smolyak1963} and developed by Barthelmann, Novak and Ritter \cite{Barthelmann2000}. The reader is referred to the extensive literature (e.g., see \cite{Novak1999,Gerstner1998,Gerstner2003}) for more details. Sparse grid stochastic collocation has been used extensively for solving parametric elliptic and parabolic PDEs \cite{Babuska2007,Nobile2008a,Xiu2005}. Here, we give only a brief outline to set up notation.
	
	First, we require a set of one-dimensional interpolation operators. Let $\{\collocSet{\level}\}_{\level=1}^{\infty}$  be a sequence of sets of $\nColloc{\level}$ distinct collocation points with $\collocSet{\level}=\{ \collocPt_\iColloc \}_{\iColloc=1}^{\nColloc{\level}} \subset \Gamma_1=[-1,1]$. For each `level' number $\level \in \mathbb{N}$ and for continuous functions $f:\Gamma_1 \to \R$, we define an interpolation operator $\interpOp{\level}:\cSpace{0}(\Gamma_1) \to  \polySpace{\nColloc{\level}-1}(\Gamma_1)$ in the usual way by
	\begin{equation}
	\interpOp{\level}[f](y):= \sum_{{\collocPt} \in \collocSet{\level}} f(\collocPt) \LagrangePoly{\collocPt}{\level}(y),
	\label{eq:interp-op-1d}
	\end{equation}
	where $\LagrangePoly{\collocPt}{\level} \in \polySpace{\nColloc{\level} -1}(\Gamma_1)$ is the Lagrange polynomial of degree $\nColloc{\level}-1$ that takes the value one at $y=\collocPt$ and zero at all the other points in $\collocSet{\level}$. If we set $\detailOp{1} := \interpOp{1}$ and define the \emph{detail operator}  $\detailOp{\level}:\cSpace{0}(\Gamma_1) \to \polySpace{\nColloc{\level}-1}(\Gamma_1)$ for $\level>1$ as
	\begin{equation}
	\detailOp{\level}:= \interpOp{\level} - \interpOp{\level-1}
	\end{equation}
	then clearly $\interpOp{\level} = \sum_{\iMi=1}^{\level} \detailOp{\iMi}$.
	
	
	Next, switching to $d$ dimensions, for any multi-index $\vMi \in \N^\nY$ we can define the so-called \emph{hierarchical surplus} (or multi-dimensional detail operator) 
	\vspace{-2mm}
	\begin{equation}
	\detailOp{\vMi} := \bigotimes_{\iY=1}^{\nY} \detailOp{\mi_\iY}.
	\end{equation}
	Given an \emph{admissible} multi-index set $\miset \subset \Nplus^\nY$, that is, a set $\miset$ such that for every $\vMi \in \miset$, 
	\begin{equation}
	\vMi - \vUnit_{j} \in \miset \text{ for } j=1,...,\nY \text{ where } \mi_j >1 \end{equation} 
	and $\vUnit_{j}:=[0,...,\underbrace{1}_{j\text{-th}},...,0]^\top,$ we define the sparse grid operator $\interpOp{\miset}$ as follows
	
	\begin{equation}
	\interpOp{\miset}[f](\y):=\sum_{\vMi \in \miset} \detailOp{\vMi}[f](\y).
	\label{eq:interp-sparse-detail}
	\end{equation}
	This is an interpolant if the sets of one-dimensional points are nested, i.e., $\collocSet{1} \subset \collocSet{2} \subset \collocSet{3} ...$ \cite[Proposition 6]{Barthelmann2000}. In that case, we can rewrite \eqref{eq:interp-sparse-detail} as 
	\begin{equation}
	\interpOp{\miset}[f](\y):=\sum_{\vCollocPt \in \collocSet{\miset}} f(\vCollocPt) \LagrangePoly{\vCollocPt}{\miset}(\y)
	\label{eq:interp-sparse-poly}
	\end{equation}
	for appropriate interpolation polynomials $\LagrangePoly{\vCollocPt}{\miset}$ {\color{otherblue} and interpolation points $\collocSet{\miset} \subset \Gamma$}. We will use Clenshaw--Curtis points.
	
	
	\begin{definition}[Clenshaw--Curtis points \cite{Clenshaw1960}]
		On the interval $\Gamma_1=[-1,1]$ the set of $\nColloc{}$ Clenshaw--Curtis (CC) points is defined by
		\begin{equation}
		\collocPt_{\iColloc} = -\cos\left(\pi\frac{\iColloc-1}{\nColloc{}-1}\right), \qquad \iColloc=1,...,\nColloc{},
		\end{equation}
		for $\nColloc{}>1$ and we set $\collocPt_{1}=0$ for $\nColloc{}=1$.
	\end{definition}
	To construct nested sets of $\nColloc{\alpha}$ CC points we set $\nColloc{1} = 1$ and choose
	\begin{equation}
	\nColloc{\level} = 2^{\level-1}+1 \text{ for } \level > 1.
	\label{eq:cc-doubling}
	\end{equation}
	
	\begin{remark}[Use of sparse grids]
		The motivation for using sparse grids is to reduce the growth in computational cost as the number $d$ of parameters increases. For tensor product grids the number of points required to achieve a fixed total degree of polynomial exactness grows exponentially. The analogous asymptotic result for sparse grids shows an optimal polynomial growth \cite[Remark 5]{Barthelmann2000}.
	\end{remark}
	
	
	\vspace{-1mm}
	
	Applying the sparse grid interpolation operator to $u^{\h}$ in \eqref{eq:full-functional-approximation} gives
	\begin{equation}
	u^{\h,\miset}(\x,t,\y) := \interpOp{\miset}[u^{\h}(\x,t,\y)].
	\label{eq:u-h-miset}
	\end{equation}
	Re-writing \eqref{eq:u-h-miset} as a linear combination of Lagrange polynomials as in \eqref{eq:interp-sparse-poly} reveals that we need to solve the ODE system \eqref{eq:system-of-odes} for each of the chosen collocation points $\vCollocPt \in \collocSet{\miset}$. The solution to each of these deterministic ODE systems needs to be further approximated using an appropriate {timestepping} algorithm. If an \emph{adaptive} scheme is used, then for each $\vCollocPt \in \collocSet{\miset}$ the timestepping method will select a (potentially different) set of timesteps $\{\timestep_{\iTime}\}_{\iTime=0}^{\nTime{}-1}$ and generate a sequence of corresponding approximations 
	\begin{equation}\vec{u}^{\h,\letol}_{\iTime}(\vCollocPt) \approx \vec{u}^{\h}(t_\iTime; \vCollocPt)
	\end{equation} 
	to the ODE solution for times $\{t_\iTime\}_{\iTime=1}^{\nTime{}} = \{\sum_{i=0}^{\iTime-1} \timestep_{i} \}_{\iTime=1}^{\nTime{}}$ where $t_{\nTime{}} \geq T$.  We will consider {\color{otherblue}stable, convergent, implicit methods} with local error control. Local error control means that we select the timestep $\timestep_\iTime$ so that the local error $\localerror{\iTime}(\vCollocPt)$ satisfies
	\begin{equation}
	\localerror{\iTime}(\vCollocPt) := \Vert \vec{u}^{\h}_{local}(t_\iTime + \timestep_\iTime, \vCollocPt) - \vec{u}^{\h,\letol}_{\iTime+1}(\vCollocPt) \Vert_{\massMatrix} \approx \delta,
	\label{eq:local-error}
	\end{equation}
	where $\letol$ is a user-chosen tolerance, $\vec{u}^{\h}_{local}(t,  \vCollocPt)$ is the ODE solution on $[t_{k}, t_{k+1}]$ with initial data  $\vec{u}^{\h}_{local}(t_\iTime, \vCollocPt) := \vec{u}^{\h,\letol}_{\iTime}(\vCollocPt)$ and $\Vert \cdot \Vert_{\massMatrix}$ is the discrete analogue of $\Vert \cdot \Vert_{\normXFull}.$
	
	After applying a timestepping method for each collocation point $\vCollocPt \in \collocSet{\miset}$, a space-time approximation ${u}^{\h,\letol}(\x,t;\vCollocPt)$ can be defined in a similar way to \eqref{eq:full-functional-approximation} by first constructing continuous-in-time coefficients from the approximations $\vec{u}^{\h,\letol}_{\iTime}(\vCollocPt)$ using linear interpolation on each subinterval $[t_\iTime, t_{\iTime+1})$ as follows
	\begin{equation}
	\vec{u}^{\h,\letol}(t;\vCollocPt) := \vec{u}^{\h,\letol}_\iTime(\vCollocPt) + \frac{t- t_\iTime}{t_{\iTime+1} - t_\iTime} \left(\vec{u}^{\h,\letol}_{\iTime+1}(\vCollocPt) - \vec{u}^{\h,\letol}_\iTime(\vCollocPt) \right).
	\label{eq:cts-in-time-at-colloc-pt}
	\end{equation}
	The final approximation to the solution of the parametric PDE \eqref{eq:adv-diff} is then
	\begin{equation}
	u^{\h,\miset,\letol}(\x,t,\y) := \sum_{\vCollocPt \in \collocSet{\miset}} u^{\h,\letol}(\x,t;\vCollocPt) \LagrangePoly{\vCollocPt}{\miset}(\y).
	\label{eq:full-approximation-scheme}
	\end{equation}
	The simple example {in \cref{sec:motivational-example}} motivates the need to adapt the multi-index set $\collocSet{\miset}$ in time as well as use adaptive timestepping.
	

	\section{Motivating Example}
	\label{sec:motivational-example}

\begin{figure}
		\begin{subfigure}[t]{0.45\textwidth}
			\centering
			\includegraphics{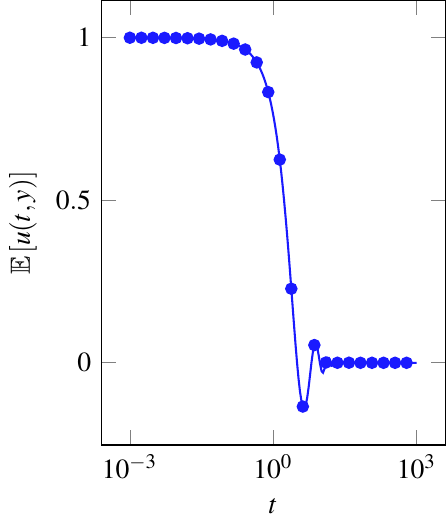}
			\caption{mean of $u(t,y)$}
			\label{fig:ode-exp}
		\end{subfigure}
		\begin{subfigure}[t]{0.45\textwidth}
			\centering
			\includegraphics{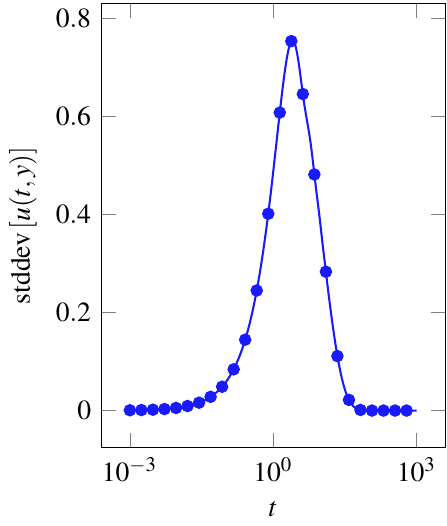}
			\caption{standard deviation of $u(t,y)$}
			\label{fig:ode-var}
		\end{subfigure}
		\caption{Mean and standard deviation of \eqref{eq:ode-explicit-sln} for $u_0=1$ and $\diffcnst=0.1$}
	\end{figure}

	Applying the Fourier transform in space to the advection-diffusion problem on $\spatialdomain=\R$ gives a system of ODEs with complex coefficients that depend on the wavenumber. We consider a single ODE of this form with a parameter-dependent coefficient $\alpha(y)$ where $y$ is the image of a uniform random variable on $\Gamma=[-1,1]$ with probability density $\rho(y)=\frac{1}{2}$.  That is, we wish to approximate $u:[0,T] \times \Gamma \to \C$ satisfying
	\begin{equation}
	\begin{cases}
	\ddt u(t,y) = \alpha(y) u(t,y) & t\in(0,T]\\
	u(0,y) = u_0\\
	\end{cases}
	\label{eq:1d-method-of-lines-ode}
	\end{equation}
	$\rho$-almost surely on $\Gamma$ where $\alpha:\Gamma \to \C$ is defined by
	\begin{equation}
	\alpha(y) = (-\diffcnst + \im y).
	\end{equation}
	Here, $\Re[\alpha(y)]=-\diffcnst$ is associated with diffusion and damps the solution in time and $\Im[\alpha(y)]=y$ represents an uncertain advective term. The exact solution to \eqref{eq:1d-method-of-lines-ode} is
	\begin{equation}
	u(t,y) = u_0 \exp(-\diffcnst t) \exp(\im y t)
	\label{eq:ode-explicit-sln}
	\end{equation}
	with mean
	\begin{equation}
	\Exp{u(t,y)} = u_0 \exp(-\diffcnst t) \sinc(t)
	\end{equation}
	and standard deviation
	\begin{equation}
	\StdDev{u(t,y)} = u_0 \exp(-\diffcnst t)\sqrt{ 1-  \sinc(t)^2}.
	\end{equation}
	These functions are plotted in \cref{fig:ode-exp,fig:ode-var} for the case $u_0=1$ and $\diffcnst=0.1$. The mean decays in time with modulation due to the $\sinc$ function. More interestingly, the standard deviation \emph{grows} in time initially.
	The solution \eqref{eq:ode-explicit-sln} becomes more oscillatory as a function of the parameter $y$ as time passes and hence will become harder to approximate in a fixed polynomial approximation space. In the long time limit the exponential decay dominates and the standard deviation decays to zero. We now apply the polynomial-based approximation scheme outlined in \cref{sec:problem-statement} and highlight the issue of long time integration for parametric problems. Issues associated with polynomial chaos approximation for a similar ODE are discussed in \cite{Wan2006}.

	
	\begin{figure}
		\begin{subfigure}[t]{0.45\textwidth}
			\centering
			\includegraphics{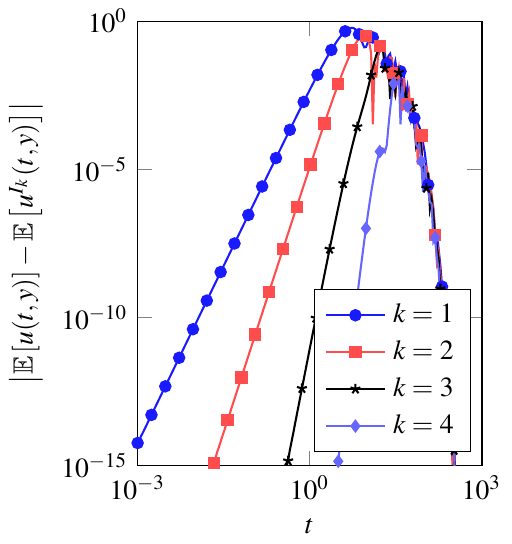}
			\subcaption{error in mean}
			\label{fig:ode-exp-error}
		\end{subfigure}
		\begin{subfigure}[t]{0.45\textwidth}
			\centering
			\includegraphics{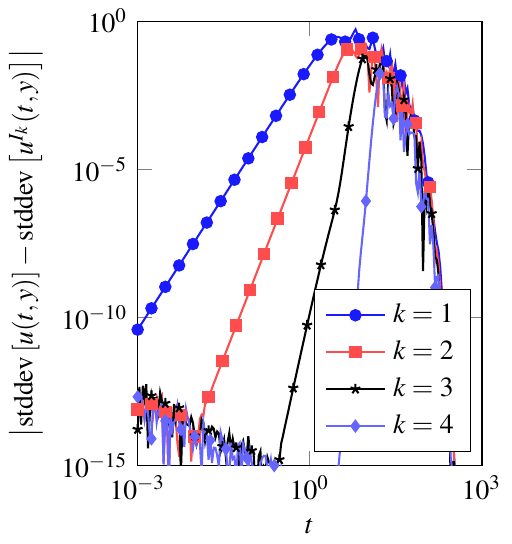}
			\caption{error in standard deviation}
			\label{fig:ode-var-error}
		\end{subfigure}
		\caption{Error in mean and standard deviation of \eqref{eq:ode-explicit-sln} for $u_0=1$ and $\diffcnst=0.1$ and increasing polynomial degree $2^\degree$}
	\end{figure}
	
	\subsection{Interpolation}
	\label{example:ode-no-timestepping}
\begin{figure}
		\includegraphics{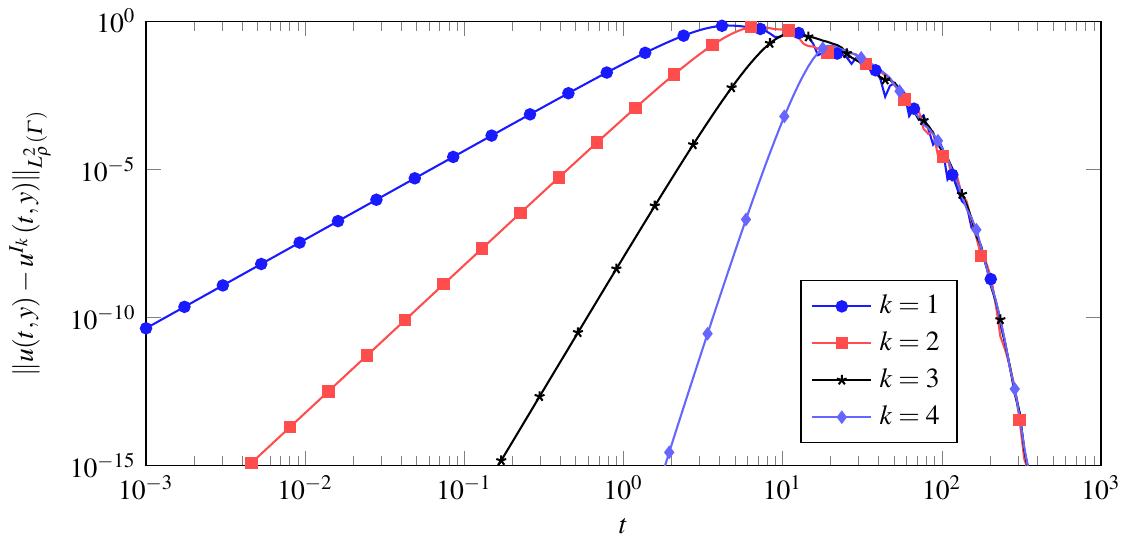}
		\caption{Approximation error $\Vert u(t,y) - u^{\miset_{\degree}}(t,y) \Vert_{\normYFull}$ for the solution of \eqref{eq:1d-method-of-lines-ode} with $u_0=1$ and $\diffcnst=0.1$ and increasing polynomial degree $2^\degree$}
		\label{fig:ode-error}
	\end{figure}

	Since we have a one-dimensional parameter domain, we define the sequence of multi-index sets $\miset_{\degree} := \{\mi\in \Nplus \st \mi \leq 1+\degree  \}$ for $\degree =1,2, \ldots$. We construct interpolants using nested sets of CC collocation points with the doubling relation \eqref{eq:cc-doubling} so that for a fixed $k$ the resulting polynomial is of degree $2^\degree$. The exact solution is available at each collocation point so timestepping is not necessary and here we simply denote approximations by $u^{\miset_{\degree}} := \interpOp{\miset_\degree}[u]$, for  $k=1,2, \ldots.$ Errors in the approximation of the mean and the standard deviation are plotted in \cref{fig:ode-exp-error,fig:ode-var-error} for $k=1,2,3,4$. These plots illustrate two important points. Firstly, for a \emph{fixed} approximation space (fixed $\degree$) the error grows with time (at least initially). Secondly, if we want to suppress the error we must grow the multi-index set $\miset_{\degree}$ in time, increasing the degree of the polynomial approximation. In \cref{fig:ode-error} the {\color{otherblue} approximation error $\Vert u(t,y) - u^{\miset_{\degree}}(t,y) \Vert_{\normYFull}$} is estimated using the $n_{11}$-point CC quadrature rule. Again the error is seen to grow in time for a fixed polynomial approximation space and can be suppressed by increasing the degree of the polynomial approximation. It is clear that adapting the set $\miset_{\degree}$ in time would be advantageous.
	


	\subsection{Timestepping}
	\label{example:ode-timestepping}
	
	\begin{figure}
		\begin{subfigure}[t]{0.45\textwidth}
			\centering
			\includegraphics{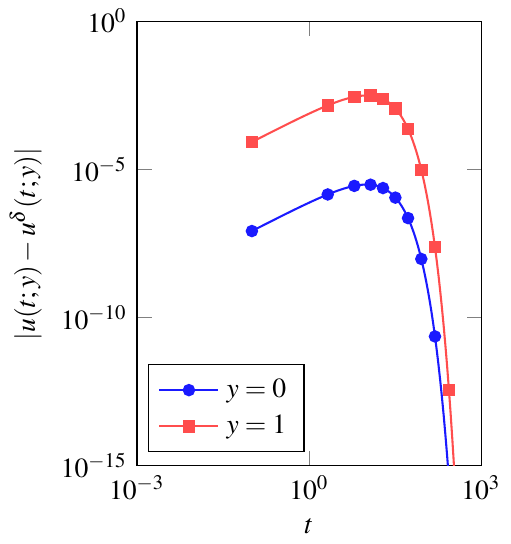}
			\caption{TR with $\timestep=10^{-1}$}
			\label{fig:ode-global-error-tr}
		\end{subfigure}
		\begin{subfigure}[t]{0.45\textwidth}
			\centering
			\includegraphics{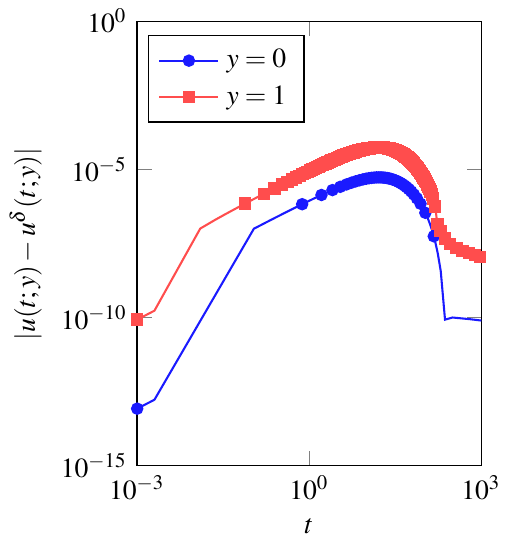}
			\caption{TR-AB2 with $\letol=10^{-7}$}
			\label{fig:ode-global-error-stabtr}
		\end{subfigure}
		\caption{Global timestepping error for approximation of \eqref{eq:1d-method-of-lines-ode} with $y=0,1$}
		\label{fig:ode-global-error}
	\end{figure}
	\begin{figure}
		\begin{subfigure}[t]{0.45\textwidth}
			\centering
			\includegraphics{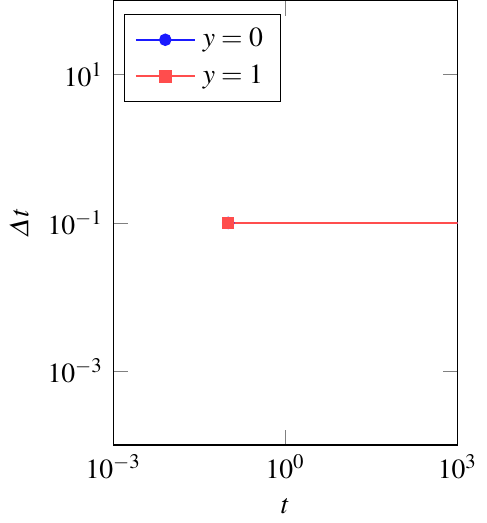}
			\caption{TR with $\timestep=10^{-1}$}
			\label{fig:ode-ts-tr}
		\end{subfigure}
		\begin{subfigure}[t]{0.45\textwidth}
			\centering
			\includegraphics{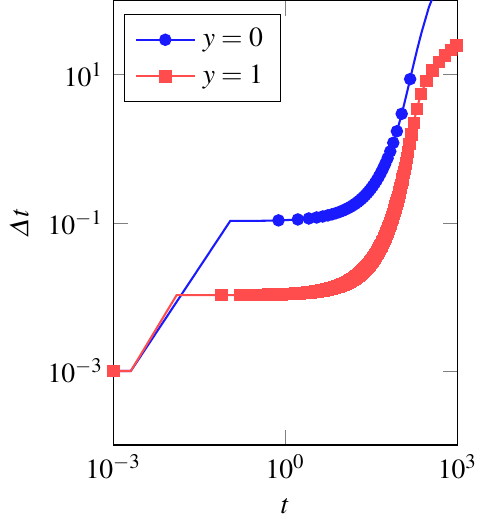}
			\caption{TR-AB2 with $\letol=10^{-7}$}
			\label{fig:ode-ts-stabtr}
		\end{subfigure}
		\caption{Timesteps for approximation of \eqref{eq:1d-method-of-lines-ode} with $y=0,1$}
		\label{fig:ode-ts}
	\end{figure}
	
	The ODE test problem \eqref{eq:1d-method-of-lines-ode} is now solved numerically using the standard {trapezoidal rule (TR)} and the \stabtr{} implementation of the {TR-AB2} algorithm \cite[Chapter 5]{Iserles2008}, \cite{Gresho2008}. {TR-AB2 pairs TR with the two-step Adams--Bashforth (AB2) method to estimate and control the local timestepping error.} We consider the deterministic ODEs associated with the fixed values $y=0$ (corresponding to $\alpha(0)=\Exp{\alpha(y)}$) and $y=1$ (which gives the maximal imaginary component of $\alpha$) and investigate the  \emph{global} timestepping error
	\begin{equation}
	\lvert u(t;y) - u^{\letol}(t;y) \rvert
	\end{equation}
	associated with the chosen schemes. The TR algorithm is applied using a fixed timestep of $\timestep = 10^{-1}$. 
	For TR-AB2, we choose a local error tolerance $\letol=10^{-7}$. The plots in \cref{fig:ode-global-error} demonstrate that for $y=0$ the two algorithms give comparable global errors for times up to $\bigO{10^2}$. In the long time limit the TR-AB2 error stagnates --- but it is already very small by this time. To reach the final time $T=10^3$, the TR method required $10^4$ steps but TR-AB2 required only $296$. This computational saving motivates the need for adaptive timestepping; we attain a comparable global error at a significantly reduced cost. \Cref{fig:ode-ts} illustrates why this saving is possible. It can be seen that in the long time limit the timestep with TR-AB2 is allowed to grow much larger than $10^{-1}$. Similar results are seen for $y=1$. In fact, a notably better global error is obtained using TR-AB2 with $\letol=10^{-7}$, requiring only $2907$ timesteps, compared to the $10^4$ timesteps used in standard TR with $\timestep=10^{-1}$.
	
	Returning to \eqref{eq:full-approximation-scheme}, we now consider how to adapt the approximation in time to ensure that the interpolation and timestepping errors remain balanced.  Key to this is an efficient a posteriori error estimator. We discuss this next.
	

\section{Error Estimation} \label{sec:error-est}
Error estimation will be performed in the norm $\Vert \cdot \Vert_{\normXY}$ where
\begin{equation}
\begin{aligned}
	\Vert u(\x,t,\y) \Vert_{\normXY} := 
	\left( \int_{\Gamma} \int_{\spatialdomain} \left(u(\x,t,\y)\right)^2 \dd \x ~\rho(\y) \dd \y \right)^\frac{1}{2}.
	\end{aligned}
\end{equation}
The $\normXFull$ norm in space is chosen as it is the most natural choice in timestepping algorithms for local error control. The $\normYFull$ norm is chosen as we consider problems with uncertain inputs and solutions that have finite variance.

To simplify the analysis the \emph{spatial} contribution to the error is not estimated --- we are interested in approximating the solution to the parametric ODE system \eqref{eq:system-of-odes}. The focus is on estimating and balancing the interpolation and timestepping errors so we consider the splitting of the error
\begin{equation}
\begin{aligned}
\Vert u^{\h} - u^{\h,\letol,\miset} \Vert_{\normXY} &\leq {\Vert u^{\h} - u^{\h,\miset} \Vert_{\normXY}} + {\Vert u^{\h,\miset} - u^{\h,\miset,\letol} \Vert_{\normXY}}\\ 
&=:{\error{\interp}(t)} +\error{\timestepping}(t).
\end{aligned}
\label{eq:error-split}
\end{equation} 
The interpolation error $\error{\interp}(t)$ is estimated using a hierarchical approach by comparing the approximation $u^{\h,\miset}$ to an \emph{enhanced approximation} and appealing to a saturation assumption. This approach follows the adaptive algorithm first proposed by Gerstner and Griebel \cite{Gerstner2003}.
To define an enhanced approximation, we first define two additional multi-index sets.
The \emph{margin} $\margin{\miset}$ of $\miset$ is defined as the set of all neighbouring multi-indices,
\begin{equation}
\margin{\miset} := \bigcup_{\vMi \in \miset} \{\vMi + \vUnit_{j} \text{ for } j\in 1,...,\nY\}.
\end{equation}
The \emph{reduced margin} $\reducedmargin{\miset}$ is defined as the set of multi-indices in $\margin{\miset}$ that can be individually added to $\miset$ whilst maintaining the admissibility property. That is, the set
\begin{equation}
\reducedmargin{\miset}:= \{\vMi \in \margin{\miset} \st \miset \cup \vMi \text{ is admissible}\}.
\end{equation}
The enhanced multi-index set is defined as $\enhanced{\miset}:= \miset \cup \reducedmargin{\miset}$ and is used to construct the enhanced approximation $u^{\h,\enhanced{\miset}}$. 
The error estimation strategy we adopt is underpinned by the following assumption.
\begin{assumption}[Saturation Assumption] \label{assumption:saturation}
	Let $\enhanced{\miset}:= \miset \cup \reducedmargin{\miset}$ be the enhanced multi-index set. There exists a saturation constant $\saturation(t) \in [0,1)$ for each $t\in[0,T]$ such that{
	\begin{equation}
	\begin{aligned}
		\Vert u^{\h}({\cdot},t,{\cdot}) - & u^{\h,\enhanced{\miset}}({\cdot},t,{\cdot}) \Vert_{\normXY} \\ &\leq \saturation(t) \Vert u^{\h}({\cdot},t,{\cdot}) - u^{{\h,\miset}}({\cdot},t,{\cdot}) \Vert_{\normXY}. \label{eq:assumption}
	\end{aligned}
	\end{equation}}
\end{assumption}

{A saturation assumption for the interpolation error estimation is a standard ingredient in a posteriori error analysis of elliptic PDEs with random data, see e.g.\,\cite[Section~4]{BSX22}. The assumption can be justified by considering the analyticity of the solution with respect to ${\bf y}$, see e.g.\,\cite{Feischl2021}.}

\begin{corollary}[Interpolation Error Bound] \label{cor:error-bound}
	If \cref{assumption:saturation} holds, then the following upper bound holds for the interpolation error,{
	\begin{equation}
	\begin{aligned}
		\error{\interp}(t) &\leq \frac{1}{1 - \saturation(t)} \Vert u^{\h,\enhanced{\miset}}(\cdot,t,\cdot) - u^{\h,{\miset}}(\cdot,t,\cdot) \Vert_{\normXY}.
	\end{aligned}
	\label{eq:interp-upper-bound}
	\end{equation}}
\end{corollary}
{
\begin{proof}
	The result follows by writing
	\begin{equation}
		\begin{aligned}
		\Vert u^{\h}(\cdot,t,\cdot) -  u^{\h,{\miset}}(\cdot,t,\cdot) \Vert_{\normXY} =&\\
		\Vert u^{\h}(\cdot,t,\cdot) -  u^{\h,\enhanced{\miset}}(\cdot,t,\cdot) +& u^{\h,\enhanced{\miset}}(\cdot,t,\cdot) - u^{\h,{\miset}}(\cdot,t,\cdot) \Vert_{\normXY}
		\end{aligned}
	\end{equation}
	and applying the triangle inequality and \cref{assumption:saturation}.
\end{proof}}

The upper bound in \eqref{eq:interp-upper-bound} is not computable --- the saturation constant is unknown and the approximations $u^{\h,\miset}(\x,t,\y)$ and $u^{\h,\enhanced{\miset}}(\x,t,\y)$ are only known \emph{after} further approximation with timestepping. Hence, we consider the further splitting
\begin{equation}
{
\begin{aligned}
\Vert u^{\h,\enhanced{\miset}}(\cdot&,t,\cdot) - u^{\h,\miset}(\cdot,t,\cdot) \Vert_{\normXY}\\
\leq &\underbrace{\Vert u^{\h, \enhanced{\miset},\letol}(\cdot,t,\cdot) - u^{\h,\miset,\letol}(\cdot,t,\cdot) \Vert_{\normXY}}_{\errorest{\miset}(t)}\\
& + \underbrace{\Vert u^{\h, \enhanced{\miset}}(\cdot,t,\cdot) - u^{\h, \enhanced{\miset},\letol}(\cdot,t,\cdot)  - u^{\h,\miset}(\cdot,t,\cdot) + u^{\h,\miset,\letol}(\cdot,t,\cdot) \Vert_{\normXY}}_{correction}
\end{aligned}
\label{eq:interp-est-and-correction}}
\end{equation}
where $\errorest{\miset}:[0,T]\to\R$ is computable but there is now an additional correction term to estimate. We will call $\errorest{\miset}{(t)}$ the interpolation error estimate.

\begin{remark}[Interpolation error indicators]
	We can bound the interpolation error estimate $\pi_{\miset}$ as
	\begin{equation}
	\errorest{\miset}(t) \leq \sum_{\vMi \in \reducedmargin{\miset}} \errorest{\miset,\vMi}(t) 
	\end{equation}
	where
	\begin{equation}
	\begin{aligned}
	\errorest{\miset,\vMi}(t):= \bigg\Vert\sum_{{\collocPt} \in \collocSet{{\miset \cup \vMi}}} u^{\h,\letol}(\x,t;\vCollocPt) \LagrangePoly{\vCollocPt}{\miset \cup \vMi}(\y)
	- \sum_{{\collocPt} \in \collocSet{{\miset}}} u^{\h,\letol}(\x,t;\vCollocPt) \LagrangePoly{\vCollocPt}{\miset}(\y)\bigg\Vert_{\normXY}.
	\end{aligned}
	\label{eq:pi-mi-alpha}
	\end{equation}
	These indicators will guide the adaptive algorithm introduced in \cref{sec:error-control}.
\end{remark}

Now consider $u^{\h,{\miset}},u^{\h,\miset,\letol},u^{\h,\enhanced{\miset}}$ and $u^{\h,\enhanced{\miset},\letol}$ in the correction term in \eqref{eq:interp-est-and-correction}. Rewriting these terms as linear combinations of Lagrange polynomials as in \eqref{eq:interp-sparse-poly} gives
\begin{equation}
\begin{aligned}
& \big(u^{\h, \enhanced{\miset}}(\x,t,\y) - u^{\h, \enhanced{\miset},\letol}(\x,t,\y)\big) - \big( u^{\h,\miset}(\x,t,\y) - u^{\h,\miset,\letol}(\x,t,\y) \big)
\\
& =\sum_{{\collocPt} \in \collocSet{\enhanced{\miset}} \setminus \collocSet{{\miset}}} e^{\h,\letol}(\x,t;\vCollocPt) \LagrangePoly{\vCollocPt}{\enhanced{\miset}}(\y) + \sum_{\collocPt \in \collocSet{{\miset}}} e^{\h,\letol}(\x,t;\vCollocPt) \left(\LagrangePoly{\vCollocPt}{{\enhanced{\miset}}}(\y) - \LagrangePoly{\vCollocPt}{{\miset}}(\y)\right),
\end{aligned}
\label{eq:correction}
\end{equation}
where we have now introduced the \emph{global timestepping error} \begin{equation}
e^{\h,\letol}(\x,t;\vCollocPt) := u^{\h}(\x,t;\vCollocPt) - u^{\h,\letol}(\x,t;\vCollocPt)
\end{equation} for each collocation point $\vCollocPt$ at all times $t\in[0,T]$. Note that this error represents both timestepping error and an interpolation error due to the extension of the approximation from discrete times $\{t_\iTime\}_{\iTime=0}^{\nTime{}}$ to $[0,t_{\nTime{}}]$ in \eqref{eq:cts-in-time-at-colloc-pt}. Let $\error{\gerror}(t;\vCollocPt):=\Vert e^{\h,\letol}(\x,t;\vCollocPt) \Vert_{\normX}$.
An upper bound for the correction term is then given by
\begin{equation}
\begin{aligned}
\Vert u^{\h, \enhanced{\miset}}(\x,t,\y) &- u^{\h, \enhanced{\miset},\letol}(\x,t,\y) + u^{\h,\miset}(\x,t,\y) - u^{\h,\miset,\letol}(\x,t,\y) \Vert_{\normXY} \\
\leq & \sum_{{\collocPt} \in \collocSet{\enhanced{\miset}} \setminus \collocSet{{\miset}}} \error{\gerror}(t;\vCollocPt) \Vert \LagrangePoly{\vCollocPt}{\enhanced{\miset}}(\y)\Vert_{\normY} \\
&+ \sum_{{\collocPt} \in \collocSet{{\miset}}} \error{\gerror}(t;\vCollocPt) \Vert \LagrangePoly{\vCollocPt}{\enhanced{\miset}}(\y) - \LagrangePoly{\vCollocPt}{{\miset}}(\y)\Vert_{\normY}\\
&=:\errorbound{\correction}(t).
\end{aligned}
\label{eq:correction-bound}
\end{equation}
A similar approach can be used to bound the timestepping error in \eqref{eq:error-split} as follows
\begin{equation}
\error{\timestepping}(t) \leq \errorbound{\timestepping}(t) := \sum_{{\collocPt} \in \collocSet{{\miset}}} \error{\gerror}(t; \vCollocPt) \Vert \LagrangePoly{\vCollocPt}{{\miset}}(\y)\Vert_{\normY}.
\label{eq:timestepping-bound}
\end{equation}
For each $\vCollocPt \in \collocSet{\enhanced{\miset}}$, let $\errorest{ge}(t; \vCollocPt):[0,T]\to\R$ be a global timestepping error estimate,
\begin{equation}
	\errorest{\gerror}(t; \vCollocPt) \approx \error{\gerror}(t;\vCollocPt). 
\end{equation}
Then, using the set of error estimates $\{\errorest{\gerror}(t;\vCollocPt)\}_{\vCollocPt \in \collocSet{\enhanced{\miset}}}$ the correction bound \eqref{eq:correction-bound} can be estimated by 
\begin{equation}
\begin{aligned}
\errorest{\miset,\letol}(t) :=\sum_{{\collocPt} \in \collocSet{\enhanced{\miset}} \setminus \collocSet{{\miset}}} & \errorest{ge}(t; \vCollocPt) \Vert \LagrangePoly{\vCollocPt}{\enhanced{\miset}}(\y)\Vert_{\normY} 
\\&+ \sum_{{\collocPt} \in \collocSet{{\miset}}} \errorest{ge}(t; \vCollocPt) \Vert \LagrangePoly{\vCollocPt}{\enhanced{\miset}}(\y) - \LagrangePoly{\vCollocPt}{{\miset}}(\y)\Vert_{\normY}
\end{aligned}
\label{eq:correction-est}
\end{equation}
and the timestepping bound \eqref{eq:timestepping-bound} can be estimated by 
\begin{equation}
\errorest{\letol}(t) := \sum_{{\collocPt} \in \collocSet{{\miset}}} \errorest{\gerror}(t; \vCollocPt) \Vert \LagrangePoly{\vCollocPt}{{\miset}}(\y)\Vert_{\normY}.
\label{eq:timestepping-est}
\end{equation}

Returning to the full bound \eqref{eq:error-split} we can write
\begin{equation}
	\Vert u^{\h} - u^{\h,\miset,\letol} \Vert_{\normXY} \leq \frac{1}{1-\saturation(t)}\left(\errorest{\miset}(t) + \errorbound{\correction}(t)\right) + \errorbound{\timestepping}(t).
	\label{eq:bound-with-saturation}
\end{equation}
Dropping the unknown $(1-\saturation(t))^{-1}$ and replacing $\errorbound{\correction}(t)$ and $\errorbound{\timestepping}(t)$ with the computable estimates $\errorest{\correction}(t)$ and $\errorest{\timestepping}(t)$, we arrive at the final error estimator
\begin{equation}
{\errorest{}(t):=  \errorest{\miset}(t) + \errorest{\miset,\letol}(t)  + \errorest{\letol}(t), \quad \quad \quad t \in [0,T].}
\label{eq:total-error-est}
\end{equation}


\subsection{Global Timestepping Error Estimation} \label{remark:global-error-estimation}
	The error estimators $\errorest{\miset,\letol}$ and $\errorest{\letol}$ defined in \eqref{eq:correction-est} and \eqref{eq:timestepping-est} require a global timestepping error estimate for each collocation point $\vCollocPt\in\collocSet{\enhanced{\miset}}$. Relating the global error incurred by timestepping algorithms to the local error is a challenging problem with an array of literature on the subject \cite{Shampine1976b,Skeel1986,Cao2004}.
	A simple approach is to compare the approximation computed with local error tolerance $\delta$ to a reference approximation computed with a smaller local error tolerance $\enhanced{\letol}\ll \letol$ and estimate the global error as follows \cite{Skeel1986} 
	\begin{equation}
	\errorest{ge}(t; \vCollocPt) := \Vert u^{\h,\enhanced{\letol}}(\x,t;\vCollocPt) - u^{\h,{\letol}}(\x,t;\vCollocPt) \Vert_{\normX}.
	\end{equation}
	This is an expensive method --- the enhanced approximation is even more expensive to compute than the original one.
	The cost can be reduced by using a scaling argument which we describe next. 
	
	We assume that we are using a stable, convergent, order $\odeP \geq 1$ time integrator with local error control and that the ODE system \eqref{eq:system-of-odes} rewritten in the form \begin{equation}
		\frac{\dd \vec{u}^{\h}}{\dd t}(t;\vCollocPt) = \vec{g}\left(t,\vec{u}^{\h}(t;\vCollocPt)\right)
	\end{equation}
	satisfies the Lipschitz condition where the solution is $\odeP+1$ times continuously differentiable.
	The global error is then  $\error{\gerror}(t;\vCollocPt)=\bigO{\timestep^\odeP}$ \cite[Definition 2.3]{Butcher2008}. Adaptive timestepping methods generally select timesteps based on the relationship between the local truncation error and timestep. For an order $\odeP$ method the local truncation error is $\bigO{\timestep^{\odeP+1}}$ and for a fixed local error tolerance $\letol$ the timesteps are chosen so that $\timestep \propto \letol^\frac{1}{\odeP+1}$ \cite{Shampine2005,Gresho2008}. Hence,
	\begin{equation}
		\error{\gerror}(t;\vCollocPt) = \bigO{\letol^{\frac{\odeP}{\odeP+1}}}.
		\label{eq:big-o-ge}
	\end{equation}
	Using \eqref{eq:big-o-ge}, for two local error tolerances $\letol$ and $\letol_0 \gg \letol$ we may assume the crude relationship
	\begin{equation}
		\frac{\error{\gerror,\letol}(t;\vCollocPt)}{\error{\gerror,\letol_0}(t;\vCollocPt)} \approx \left(\frac{\letol}{\letol_0}\right)^\frac{\odeP}{\odeP+1},
	\end{equation}
	where the superscript denotes the global error corresponding to the local error tolerance $\letol$, and similarly for $\letol_0$. Therefore,
	\begin{equation}
	{\error{\gerror,\letol}(t;\vCollocPt)} \approx \left(\frac{\letol}{\letol_0}\right)^\frac{\odeP}{\odeP+1} \error{\gerror,\letol_0}(t;\vCollocPt).
	\label{eq:scaling}
	\end{equation}
	
	Using a low fidelity approximation $u^{\h,\letol_0}$ with tolerance $\letol_0 \gg \letol$ and comparing to the computed approximation $u^{\h,\letol}$ means a global error estimate for the tolerance $\letol_0$ can be formed as
	\begin{equation}
		\errorest{ge,\letol_0}(t,\vCollocPt):=\Vert u^{\h,{\letol}}(\x,t;\vCollocPt) - u^{\h,{\letol_0}}(\x,t;\vCollocPt) \Vert_{\normX} \approx \error{\gerror,\letol_0}(t,\vCollocPt).
	\end{equation}
	Combining this with \eqref{eq:scaling} yields a cheap to compute global error estimate for approximations computed with the local error tolerance $\letol$. That is,
	\begin{equation}
	\begin{aligned}
	\errorest{\gerror}(t; \vCollocPt) := & \left(\frac{\letol}{\letol_0}\right)^{\frac{\odeP}{\odeP+1}} \errorest{\gerror,\letol_0}(t,\vCollocPt)\\
	=&\left(\frac{\letol}{\letol_0}\right)^{\frac{\odeP}{\odeP+1}}\Vert u^{\h,{\letol}}(\x,t;\vCollocPt) - u^{\h,{\letol_0}}(\x,t;\vCollocPt) \Vert_{\normX}.
	\end{aligned}
	\label{eq:global-timestepping-error-estimate}
	\end{equation}
	We can then use this to compute our error estimates $\errorest{\correction}(t)$ and $\errorest{\timestepping}(t)$ in \eqref{eq:correction-est} and \eqref{eq:timestepping-est}.
	
	\begin{remark}[Choice of timestepping method]
		In our numerical experiments we will use the \stabtr{} implementation of TR-AB2 that is included in the \matlab{} software package \ifiss{} \cite{Elman2007}. This is an order $\odeP=2$ method. Many alternative algorithms with local error control exist, including the built{-}in \matlab{} integrators described in \cite{Shampine1997}. Ideally, we would use a time integration method that includes global error estimation but these are rarely implemented and used in practice. An example of such a method is \textsc{GERK} \cite{Shampine1976} but this is both an explicit method and more expensive than the TR-AB2 integrator.
	\end{remark}

\section{Adaptive Strategy} \label{sec:error-control}
To efficiently construct an approximation to the solution of the semi-discrete parametric advection--diffusion problem \eqref{eq:system-of-odes}, we will \emph{adapt} the multi-index set $\miset$ in time whilst allowing the timestepping algorithm to adaptively choose timesteps for the deterministic systems associated with individual collocation points. The construction \eqref{eq:full-approximation-scheme} only requires a small modification to allow this.
The approximation is now constructed with a {time-dependent} multi-index set \begin{equation}\miset:[0,T] \to \{ \text{admissible multi-index sets in }\nY\text{ parameters}\}\end{equation} with the adaptive approximation defined pointwise in time as
\begin{equation}
u^{\adapt}(\x,t,\y) := \sum_{{\vCollocPt} \in \collocSet{{\miset(t)}}} u^{\h,\letol}(\x,t;\vCollocPt) \LagrangePoly{\vCollocPt}{\miset(t)}(\y).
\label{eq:full-approximation-scheme-adaptive}
\end{equation}

To construct such an approximation we use a solve--estimate--mark--refine type loop in the spirit of Gerstner and Griebel to select an appropriate multi-index set \cite{Gerstner2003}.
At carefully selected times, we {estimate} the interpolation error by $\errorest{\interp}(t)$ as defined in \eqref{eq:interp-est-and-correction}, noting that the multi-index set $\miset$ can now change through time. If the estimated interpolation error is greater than a chosen tolerance $\interpTol$ at time $t$, we {mark} a set of multi-indices $\marked \subset \reducedmargin{\miset(t)}$ based on the error indicators $\{\errorest{\miset,\vMi}(t)\}_{\vMi\in\reducedmargin{\miset(t)}}$ defined in \eqref{eq:pi-mi-alpha}. The approximation \eqref{eq:full-approximation-scheme-adaptive} is then {refined} by adding the marked multi-indices to the construction of the approximation.

The \rbl{time dependence} in the problem provides an additional challenge as this adaptive strategy must be embedded 
\rbl{inside} a time propagation loop. The approximation at each collocation point is constructed independently with different discrete timesteps. To 
compute error estimates, the set of approximations must be synchronised at a common set of \emph{synchronisation times}.
\rbl{Such a set of synchronization times could be specified a priori, but would be better  generated dynamically.
	This can be done by  starting at time $t=0$ with an initial 
	synchronization timestep $\timestepAlg:=\timestepAlg_0$. 
	The approximation \eqref{eq:full-approximation-scheme-adaptive} and the error estimators $\errorest{\interp}$, 
	$\errorest{\correction}$  and $\errorest{\timestepping}$ defined in 
	\eqref{eq:interp-est-and-correction}, \eqref{eq:correction-est} and \eqref{eq:timestepping-est}  would  then be
	computed  at the new time $t+\timestepAlg$.  
	At this new time we must decide if the multi-index set needs to be refined by comparing $\errorest{\interp}(t+\timestepAlg)$ to a chosen tolerance $\interpTol(t+\timestepAlg)$. 
	Given that  the goal is to  balance the error contributions  then a
	sensible strategy is to link the interpolation error tolerance
	to the correction error estimated by \eqref{eq:correction-est}. Specifically, we define
	\begin{equation} \label{eq:toldef}
	\interpTol(t) := \scalingInterpTol \cdot \errorest{\miset,\letol}(t),
	\end{equation}
	with a \emph{safety factor} $\scalingInterpTol>1$.
	If the computed interpolation error estimate $\errorest{\miset}(t+\timestepAlg)$ 
	is   less than $\interpTol(t+\timestepAlg)$ then the synchronisation step would be  \emph{accepted}.
	In this case, the time is reset  $t \gets t+\timestepAlg$ and 
	the synchronization timestep  is increased  $\timestepAlg \gets \scalingUp\timestepAlg$,
	where $\scalingUp \geq 1$ is a user-defined input. 
	Otherwise, the step is  \emph{rejected}.  
	In this case the interpolant would be  {refined} so that  $\miset \gets \miset \cup \marked$
	and the synchronization timestep  is decreased $\timestepAlg \gets \scalingDown \timestepAlg$, 
	where $\scalingDown \in (0,1)$ is a second user-defined input. }
\rbl{
	The attraction of this  construction is that it generates an approximation in which the 
	interpolation error is {\it balanced} with the timestepping  error.  The strategy does not attempt
	to \emph{control} the global timestepping error.
}

\rbl{
	In our experiments, time integration over each interval $[t,t+\timestepAlg]$ is done  by running
	TR--AB2  adaptive timestepping with a local error tolerance $\letol$. 
	The TR--AB2 algorithm outputs a sequence of spatial approximations and the corresponding 
	intermediate timesteps for every collocation point. 
	On each interval, for each collocation point, the adaptive timestepping will be initialised with the corresponding solution approximation at the previous two time integrator steps.  This is necessary as AB2 is not self-starting.
	Global  error estimates are computed using  \eqref{eq:global-timestepping-error-estimate}
	as discussed in \cref{remark:global-error-estimation}.
}

\rbl{The complete algorithm specification will also require a marking strategy for the multi-index sets.}
A generic  strategy, similar to that used in many related algorithms \cite{Dorfler1996}, is 
\rbl{specified  in} \cref{alg:marking}. The inputs are the \rbl{current} reduced margin $\reducedmargin{\miset}$, the
corresponding error indicators $\{\errorest{\miset,\vMi_{\iMi}}\}_{\iMi=1}^{\nMi}$ and 
\rbl{a threshold parameter  $\dorfler$}. 
The algorithm returns a set of marked multi-indices $\marked\subset \reducedmargin{\miset}$ that will be 
used to refine the approximation \eqref{eq:full-approximation-scheme-adaptive}.
\begin{algorithm}
	\begin{algorithmic}[1]
		\Procedure{Mark}{$\reducedmargin{{\miset}}$, $\{\errorest{\miset,\vMi_{\iMi}}\}_{\iMi=1}^{\nMi}$, $\dorfler$}
		\State{Order $\vMi \in \reducedmargin{{\miset}}$ such that $\errorest{\miset,\vMi_1} \geq \errorest{\miset,\vMi_2} \geq \cdots \geq \errorest{\miset,\vMi_\nMi}$.}
		\State{$\eta=0, \nMarked = 0$}
		\While{$\eta < (1-\theta)$}
		\State{$\nMarked \gets \nMarked + 1$}
		\State{$\eta \gets \sum_{\iMi=1}^\nMarked \errorest{\miset,\vMi_i} / \sum_{\iMi=1}^\nMi \errorest{\miset,\vMi_i} $}
		\EndWhile
		\State\Return{$\marked = \{\vMi_\iMi\}_{\iMi=1}^{\nMarked}$}
		\EndProcedure
	\end{algorithmic}
	\caption{{ D\"{o}rfler marking}}
	\label{alg:marking}
\end{algorithm}

Alternative marking schemes can be designed, for example incorporating computational cost to give a measure of 
error per unit of computational work. This would order multi-indices by the profit indicators 
$\errorest{\miset,\vMi}/\nColloc{\vMi}$ where $\nColloc{\vMi}$ is the number of {additional} collocation points 
required; for further details see, for example,~\cite{Nobile2016,Guignard2018,Eigel2022}. 

The complete solution algorithm is given in \cref{alg:adaptive-loop}. 
The input  parameters are the local error tolerance $\letol$, the safety factor $\scalingInterpTol$, and the marking parameter  $\dorfler$.
The other parameters on \cref{ln:initialise} are the spatial discretisation parameter $\h$, the initial timestep $\timestep_0$, initial synchronisation timestep $\timestepAlg_0$, and the scaling factors $c_{\pm}$. 
They are distinguished  from the inputs. While  they affect the constructed approximation, they have little influence on  the  resulting approximation error.


{ The performance of \cref{alg:adaptive-loop} will be illustrated on two test problems in \cref{sec:numerical-results} and \cref{sub:example-sixtyfour}.}

\begin{algorithm}
	\begin{algorithmic}[1]
		\Procedure{Adaptive SC-FEM}{$\letol, \scalingInterpTol, \dorfler$}
		\State{Fix $\h,\timestep_0, \{\timestepAlg_0, \scalingDown, \scalingUp\}$.} \Comment{\rbl{Setup}} \label{ln:initialise}
		\State{Initialise multi-index sets $\miset \gets \{ \mistyle{1} \}$, $\reducedmargin{\miset}$, $\enhanced{\miset} \gets \miset \cup \reducedmargin{\miset}$.} 
		\State{Construct $\massMatrix$,$\diffMatrix$, $\advMatrix(\vCollocPt)$, $\forceFn(t,\vCollocPt)$ as defined in \eqref{eq:fem-matrices} and \eqref{eq:fem-f} for all $\vCollocPt \in \collocSet{\enhanced{\miset}}$.}
		\State{Initialise $t = 0,\timestepAlg = \timestepAlg_0$ and $u^{\h,\letol}_{local}(\x;\vCollocPt) = u^{\h}_0(\x), \timestep_\vCollocPt = \timestep_0$ for all $\vCollocPt \in \collocSet{\enhanced{\miset}}$.}
		\State{}
		\While{$t < T$}
		\For{$\vCollocPt \in \collocSet{\enhanced{\miset}}$} \Comment{Solve}
		\State{$[\{u^{h,\letol}(\x,t_\iTime; \vCollocPt)\}_{\iTime=1}^{\nTime{}}, \{\timestep_\iTime\}_{\iTime=1}^{\nTime{}}] = \textsc{Timestepping}(\letol,u^{\h,\letol}_{local}(\x;\vCollocPt), t, t + \timestepAlg, \timestep_\vCollocPt)$.}
		\State{Interpolate $u^{h,\letol}(\x,s; \vCollocPt)$ in time for $s\in[t, t+\timestepAlg]$.}
		\State{Update $u^{\h,\letol}_{local}(\x;\vCollocPt) \gets u^{h,\letol}(\x,t_{\nTime{}}; \vCollocPt)$, $\timestep_\vCollocPt \gets \timestep_{\nTime{}}$.}
		\EndFor
		\State{Construct $u^{\adapt}(\x,s,\y)=\sum_{\vCollocPt\in\miset}u^{\h,\letol}(\x,s;\vCollocPt)\LagrangePoly{\vCollocPt}{\miset}(\y)$ for $s\in[t,t+\timestepAlg]$.}
		\State{Compute $\errorest{\miset}(t+\timestepAlg),\errorest{\miset,\letol}(t+\timestepAlg),\errorest{\letol}(t+\timestepAlg)$.} \Comment{Estimate}
		\State{$\interpTol \gets \scalingInterpTol\cdot \errorest{\miset,\letol}(t+ \timestepAlg)$.}
		\State{}
		\If{$\errorest{\miset}(t+ \timestepAlg) > \interpTol(t+ \timestepAlg)$} 		\Comment{Reject}
		\State{Compute $\{ \errorest{\miset,\vMi}(t+ \timestepAlg) \}$ 	for $\vMi \in \reducedmargin{\miset}$.}
		\State{$\marked=\textsc{Mark}(\reducedmargin{\miset},\{\errorest{\miset,\vMi}(t+ \timestepAlg)\}_{\vMi \in \reducedmargin{\miset}}, \dorfler)$.} \Comment{Mark}
		\For{$\vCollocPt \in \collocSet{\enhanced{\miset} \cup \marked} \setminus \collocSet{\enhanced{\miset}}$} \Comment{Refine}
		\State{Initialise $u^{\h,\letol}_{local}(\x;\vCollocPt) = u^{\h}_0(\x), \timestep_\vCollocPt = \timestep_0$ and $\advMatrix(\vCollocPt)$, $\forceFn(t,\vCollocPt)$.} \label{ln:refine-initialise}
		\State{$[\{u^{h,\letol}(\x,t_\iTime; \vCollocPt)\}_{\iTime=1}^{\nTime{}}, \{\timestep_\iTime\}_{\iTime=1}^{\nTime{}}] = \textsc{Timestepping}(\letol,u^{\h,\letol}_{{local}}(\x,\vCollocPt),0,t,\timestep_\vCollocPt)$.} \label{ln:refine-catchup}
		\State{Update $u^{\h,\letol}_{local}(\x;\vCollocPt) \gets u^{h,\letol}(\x,t_{\nTime{}}; \vCollocPt)$, $\timestep_{\vCollocPt} \gets \timestep_{\nTime{}}$.} \label{ln:refine-update}
		\EndFor
		\State{$\miset \gets \miset \cup \marked$}
		\State{$\enhanced{\miset}\gets {\miset} \cup 	\reducedmargin{\miset}$}
		\State{$\timestepAlg \gets \scalingDown \timestepAlg$}
		\Else  \Comment{Accept}
		\State{$t \gets t+\timestepAlg$}
		\State{$\timestepAlg \gets \scalingUp\timestepAlg$}
		\EndIf
		\EndWhile
		\State{}
		\State{\Return{$u^{\adapt}(\x,t,\y)$ for $(\x,t,\y) \in \spatialdomain \times [0,T] \times \Gamma$}}
		\EndProcedure
	\end{algorithmic}
	\caption{Adaptive SC-FEM for parametric advection--diffusion}
	\label{alg:adaptive-loop}
\end{algorithm}

{
\begin{remark}[Shrinking the multi-index set]
	The diffusive term will smooth the solution through time as illustrated in the motivational example in \cref{sec:motivational-example}.
	This suggests that it may be desirable to simplify the parametric approximation by shrinking the multi-index set.
	This could be achieved by computing the error indicators \eqref{eq:pi-mi-alpha} for each $\vMi \in \miset{}$ and defining a pruning algorithm to remove the smallest contributions.
	Care would be required to ensure the resulting multi-index set is still admissible.
	In practice, the computational savings may not be significant.
	At later times when the diffusive smoothing is the dominant effect, the timesteps will have already grown large as seen in \cref{fig:ode-ts-stabtr}.
	Advancing the approximation in time, even with a large number of collocation points, is then cheap to compute.	
\end{remark}	
}

\section{{Computational {experiment}: $\nY=4$}} \label{sec:numerical-results}
{We first consider a low dimensional problem for which an accurate reference solution can be computed to test the efficiency and reliability of the error estimator in \cref{alg:adaptive-loop}.}
{We consider a generalization of the  double glazing} problem~\cite{Elman2014}.
The problem is solved on a square domain $\spatialdomain=(-1,1)^2$ with
\rbl{Dirichlet boundary data} 
\begin{equation}
u_{\spatialbdry}(\x,t) := \begin{cases} 0 & \x \in \spatialbdry_0\\
({1-x_2^4}) \left(1- \exp\left({-t}/\bdryrate\right)\right) & \x \in \spatialbdry_1
\end{cases} ,
\label{eq:hotwall}
\end{equation}
where $\spatialbdry_1 = \{\x \in \spatialbdry \st x_1 = 1\}$ \rbl{and}
$\spatialbdry_0:= \spatialbdry \setminus \spatialbdry_1 $.  
\rbl{The parameter $\bdryrate$ in \eqref{eq:hotwall} controls the rate at which the 
	right-hand boundary  wall} `heats up'. 
The initial condition is identically zero,
\begin{equation}
u(\x,0,\y)=u_0(\x,\y)=0.
\end{equation}
The mean wind field is \rbl{given by}
\begin{equation}
\wind_0\rbl{(x_1,x_2)} =  \left(2 x_2(1-x_1^2), -2 x_1(1-x_2^2)\right)
\label{eq:w0_def}
\end{equation}
and uncertainty will be introduced by perturbing $\wind_0$.



\rbl{The parameter domain is given by} $\Gamma = [-1,1]^4$ and  \rbl{we set} $\pdensity(\y) = (0.5)^4$. 
The spatial domain $\spatialdomain$ is \rbl{subdivided} into  four quadrants
\begin{equation}
\spatialdomain_\iY = \left[a_\iY-0.5, a_\iY+0.5\right]\times\left[b_\iY-0.5, b_\iY + 0.5\right] ,
\label{eq:subdomains-4}
\end{equation} where $\vec{a}=[-0.5,-0.5,0.5,0.5]^\top$, $\vec{b} = [-0.5,0.5,-0.5,0.5]^\top$.
The parametric wind field  $\wind:\spatialdomain \times \Gamma \to \R^2$   is \rbl{then defined 
	in piecewise fashion so that}
\begin{equation}
\wind(\x,\y) = \wind_0(\x) + \sigma \sum_{\iY=1}^{\rbl{4}} y_\iY \,  \wind_{\iY}(\x)
\label{eq:eddies-wind}
\end{equation}
where $\sigma = 0.5$, \rbl{with smaller {locally divergence-free} recirculations given by}
\begin{equation}
\wind_\iY(\x):= \begin{cases} \rbl{ \wind_0(z_1,z_2),} &  \rbl{z_1=2( x_1-a_\iY), z_2= 2(x_2 - b_\iY), }
\quad \x  \in \spatialdomain_ \iY\\
0, & \text{otherwise.}
\end{cases} 		
\end{equation}
\rbl{Note that, while the linearity of \eqref{eq:eddies-wind} ensures that the perturbed wind field remains divergence-free, the tangential components  of the perturbed wind are discontinuous along the interior subdomain boundaries.
By design, the non-uniform finite element {mesh} is aligned with these lines of discontinuity.}

The spatial  \rbl{discretization  is generated  by~\ifiss{} using continuous piecewise bilinear approximation.} 
This gives a parametric system of ODEs as specified in \eqref{eq:system-of-odes}.
Timestepping is performed using \stabtr{} as defined in \cite{Gresho2008}. 
Approximation input {options} are specified in \cref{tab:ifissAndStabtr}. 
Nested sets of CC collocation points \rbl{are}  generated using \sparsekit{} \cite{Back2011,Piazzola2022}. 

\begin{table}[!th]
	\caption{\ifiss{} and \stabtr{} {input options used for approximating the parametric double glazing problem with $\nY=4$}}
	\label{tab:ifissAndStabtr}
	\begin{tabular}[t]{ll}
		\hline\noalign{\smallskip}
		\ifiss{} option & value\\
		\noalign{\smallskip}\hline\noalign{\smallskip}
		Hot wall time constant $\bdryrate$ & $10^{-1}$\\
		Grid parameter $\gridParam$  ($2^{2\gridParam}$ elements)& 4\\ 
		Stretched grid & yes\\
		\rbl{Viscosity parameter $\epsilon$} & \rbl{$10^{-1}$ }\\
		\noalign{\smallskip}\hline
	\end{tabular} \qquad
	\begin{tabular}[t]{ll}
		\hline\noalign{\smallskip}
		\stabtr{} option & value\\
		\noalign{\smallskip}\hline\noalign{\smallskip}
		Averaging counter $\nstar$ & no averaging \\
		Initial timestep $\timestep_0$ & $10^{-9}$\\
		Local error tolerance $\letol$ & $10^{-5}$\\
		Final time $T$ & $100$ \\
		\noalign{\smallskip}\hline
	\end{tabular}
\end{table}
The error in the \textsc{Adaptive SC-FEM} approximation will be evaluated at a set of fifty 	logarithmically-spaced \emph{reference times} $\refTimes$ 
between $\bdryrate \log\left(0.9^{-1}\right)$ and $10^2$ and computed with respect to a \emph{reference}  approximation $u^{\reference}$ as
\begin{align}
\Vert u^{\reference}- u^{\adapt} \Vert_{\normXY} \approx \error{\adapt}: =  \Vert u^{\h}- u^{\adapt} \Vert_{\normXY}.
\end{align}
The reference approximation $u^{\reference}$ is constructed using a fixed multi-index set together with high-fidelity timestepping with the specific {input options} listed in \cref{table:ref-params}.

\begin{table}[!hb]
	\caption{Reference approximation {input options used for the parametric double glazing problem with $\nY=4$}}
	\label{table:ref-params}
	\begin{tabular}{l l}
		\hline\noalign{\smallskip}
		option & value\\
		\noalign{\smallskip}\hline\noalign{\smallskip}
		Reference multi-index set $\miset_{\reference}$ & $\{\vMi \st \Vert \vMi \Vert_1 \leq 4 + 5 \}$\\
		Reference local error tolerance $\letol_{\reference}$ & $10^{-7}$\\
		\noalign{\smallskip}\hline
	\end{tabular}
\end{table}

\begin{figure}[!t]
	\begin{subfigure}{0.49\textwidth}
		\centering
		\includegraphics{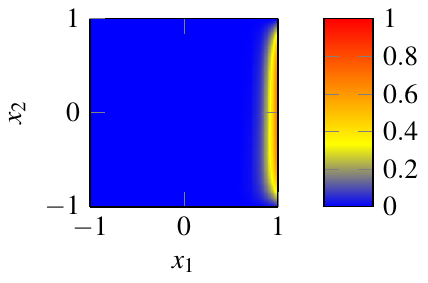}
		\caption{$\quad t\approx0.1$}
	\end{subfigure}\hspace{-10pt}
	\begin{subfigure}{0.49\textwidth}
		\centering
		\includegraphics{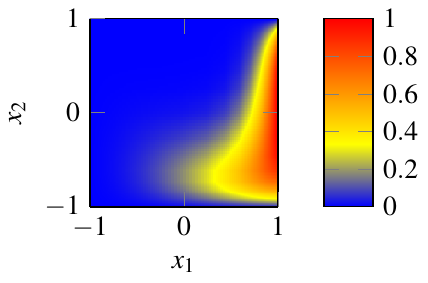}
		\caption{\quad $t\approx1$}
	\end{subfigure}
	\begin{subfigure}{0.49\textwidth}
		\centering
		\includegraphics{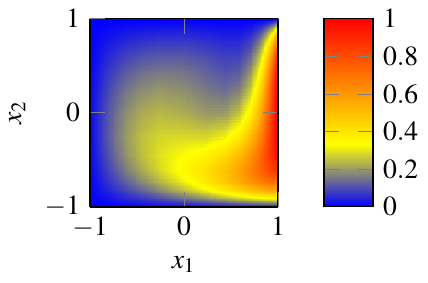}
		\caption{\quad $t\approx10$}
	\end{subfigure}\hspace{-10pt}
	\begin{subfigure}{0.49\textwidth}
		\centering
		\includegraphics{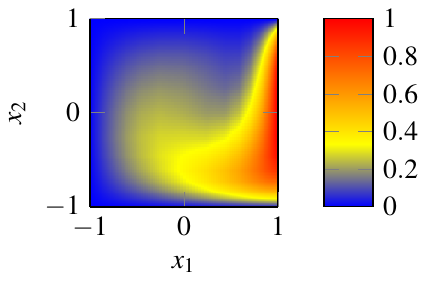}
		\caption{$\quad t\approx100$}
	\end{subfigure}	\caption{\rbl{Snapshots in time of the mean of the reference approximation when solving the parametric double glazing problem} { with $\nY=4$}}
	\label{fig:exp-pde-l4}
\end{figure}
\begin{figure}[!h]
	\begin{subfigure}{0.49\textwidth}
		\centering
		\includegraphics{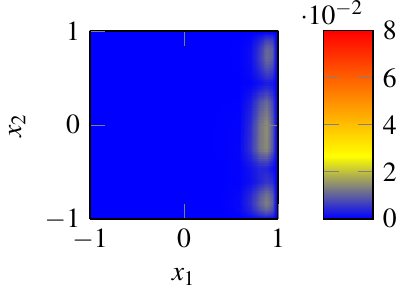}
		\caption{$\quad t\approx0.1$}
	\end{subfigure}\hspace{-10pt}
	\begin{subfigure}{0.49\textwidth}
		\centering
		\includegraphics{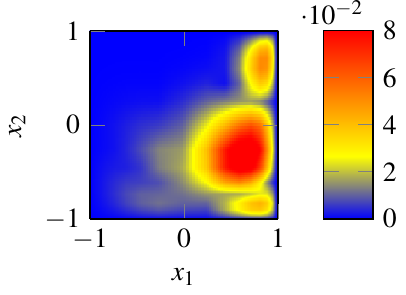}
		\caption{$\quad  t\approx1$}
	\end{subfigure}
	\begin{subfigure}{0.49\textwidth}
		\centering
		\includegraphics{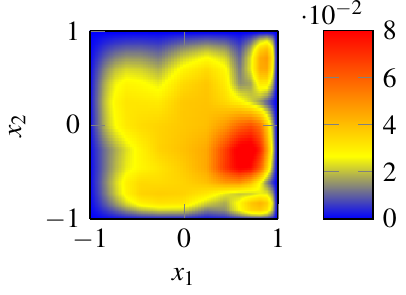}
		\caption{$\quad t\approx10$}
	\end{subfigure}\hspace{-10pt}
	\begin{subfigure}{0.49\textwidth}
		\centering
		\includegraphics{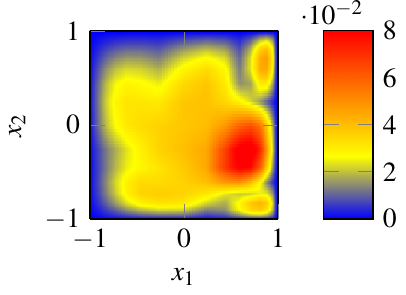}
		\caption{$\quad t\approx100$}
	\end{subfigure}
	\caption{\rbl{Snapshots in time of the standard deviation of the reference approximation when solving the parametric double glazing problem} { with $\nY=4$}} 
	\label{fig:stddev-pde-l4}
\end{figure}

\rbl{The mean of the reference approximation
	$\Exp{u^{\reference}}$  
	is shown in \cref{fig:exp-pde-l4}. The snapshots plotted} demonstrate the expected advection--diffusion behaviour: 
at short times a diffusive boundary layer forms on the ``hot wall'', as time passes heat is advected
clockwise around the domain \rbl{before eventually reaching} a steady state. 
The standard deviation  \rbl{of the reference approximation is illustrated in} \cref{fig:stddev-pde-l4}. \rbl{As for the} 
ODE problem in \cref{example:ode-no-timestepping}, \rbl{the standard deviation { at each point in space will initially} increase in magnitude as time passes: this suggests that}  the multi-index set in an adaptive parametric approximation will need to grow in time to keep the interpolation error under control.
{ In the long time limit the standard deviation reaches a steady, non-zero state.}

\begin{figure}[!ht]
	\centering
	\includegraphics{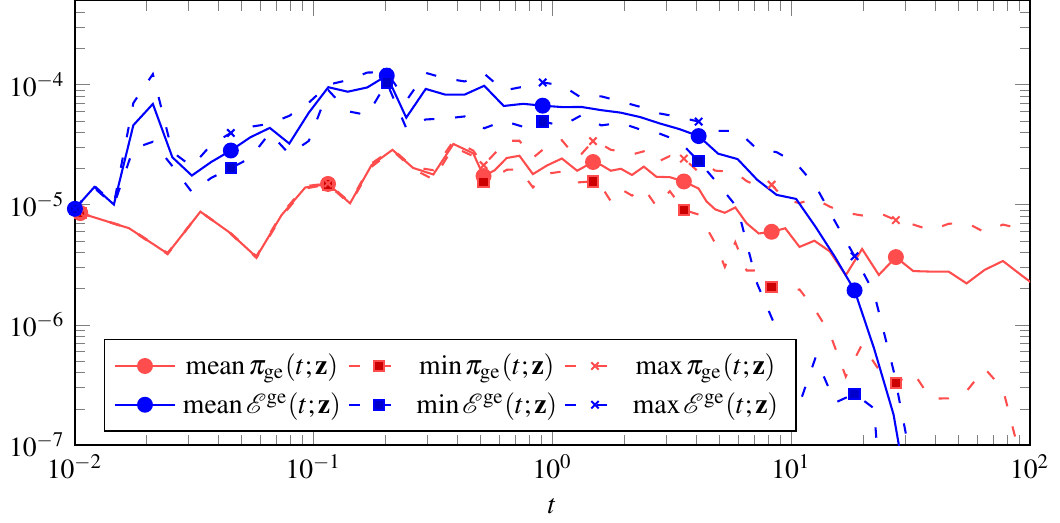}
	\caption{{Sample mean, minimum and maximum of the global timestepping error estimates $\{\errorest{\gerror}(t;\vCollocPt) \}_{\vCollocPt\in\collocSet{\adapt}}$ (in red with circle, square and cross marks, respectively) and of the reference errors $\{\error{\gerror}(t;\vCollocPt) \}_{\vCollocPt\in\collocSet{\adapt}}$ (in blue with the same markers for comparison) for the \rbl{parametric} double glazing problem  with $\nY=4$}}
	\label{fig:pde-global-error-est}
\end{figure}

\rbl{Having computed a  reference solution, we are now in a position to test
	the reliability and efficiency of the error estimation strategy in \cref{alg:adaptive-loop}. In the first instance, the}
synchronisation timestep is initialised as $\timestepAlg_0 = \bdryrate \log\left(0.9^{-1}\right)$ and  the
scaling parameters  \rbl{are fixed as} $\scalingUp=1.2$ {and} $\scalingDown=0.5$. 
Approximations are formed for collocation points with an initial timestep $\timestep_0 = 10^{-9}$. 
The \rbl{local error} tolerance is chosen to be $\letol=10^{-3}$, the safety factor $\scalingInterpTol=10^1$ 
and the marking factor $\dorfler=10^{-1}$. Different choices of these parameters will be investigated later.
First, we investigate  the global timestepping error estimate \eqref{eq:global-timestepping-error-estimate}. 
\rbl{Denoting the time-dependent set of collocation points by $\collocSet{\adapt}(t):=\collocSet{\miset{(t)}}$, the estimates $\{\errorest{\gerror}(t;\vCollocPt) \}_{\vCollocPt\in\collocSet{\adapt}}$ are computed} 
at the reference times $t\in\refTimes$.
The sample mean, minimum and maximum values of \rbl{these} estimates are plotted in \cref{fig:pde-global-error-est} . 
These are compared to the same statistics for the errors computed with respect to the reference approximation, 
$\Vert u^{\h,\letol_{\reference}} - u^{\h,\letol} \Vert_{{L}^{2}(\spatialdomain)}
\approx \error{\gerror}(t)=\Vert e^{\h,\letol} \Vert_{{L}^{2}(\spatialdomain)}$.
\rbl{Looking at these results, we see that the true error is
	generally \emph{underestimated} for short times and \emph{overestimated} when close to the steady state.} 
There is room for improvement with a more 
sophisticated global timestepping error estimator, but \eqref{eq:global-timestepping-error-estimate} will suffice for our purposes as we will see shortly, when we consider the {total} estimator, {$\errorest{}(t)$}.

\begin{figure}
	\begin{subfigure}[t]{0.49\textwidth}
		\centering
		\includegraphics{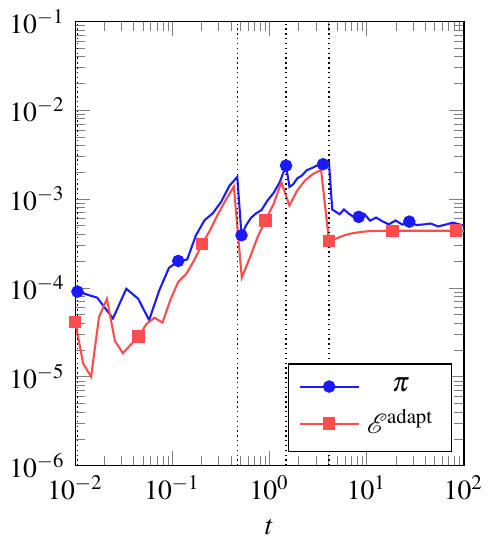}
		\caption{Error estimator $\errorest{}(t)$ and actual error $\error{\adapt}(t)$}
		\label{fig:pde-l4-error-est-nointerp}
	\end{subfigure}
	\begin{subfigure}[t]{0.49\textwidth}
		\centering
		\includegraphics{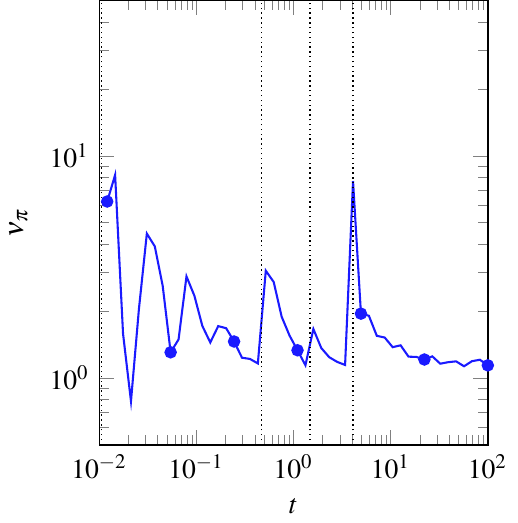}
		\caption{Efficiency $\efficiency{\errorest{}}{(t)} := {\errorest{}(t)}/{\error{\adapt}(t)}$}
		\label{fig:pde-l4-error-eff-nointerp}
	\end{subfigure}
	\caption{Error estimator $\errorest{}(t)$ for the \rbl{parametric} double glazing problem { with $\nY=4$}. 
	Vertical dotted lines denote \rbl{the parametric refinement} times}
\end{figure}

Next, we consider the quality of the \rbl{total} error estimator $\errorest{}(t)$ defined in \eqref{eq:total-error-est}.
\cref{fig:pde-l4-error-est-nointerp} demonstrates that the error estimator \rbl{provides} a close estimate of the 
true error $\error{\adapt}$. This is reflected in the \rbl{effectivity}  of the estimator 
$\efficiency{\errorest{}}(t)$ plotted in \cref{fig:pde-l4-error-eff-nointerp}. 
The error estimator generally gives a slight overestimation the true error, with the most 
notable spikes occurring at times where parametric refinement is done. 

\begin{figure}
	\centering
	\includegraphics{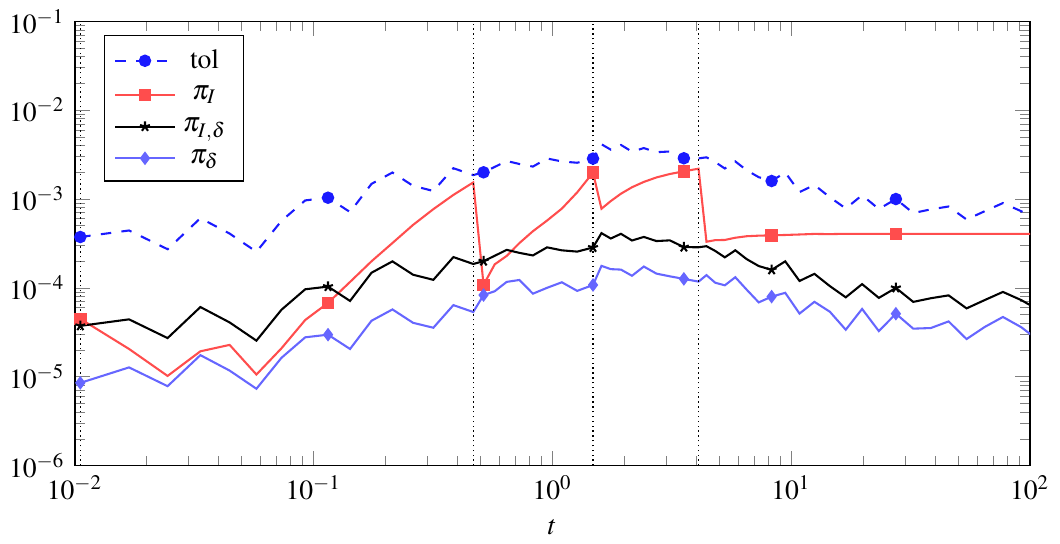}
	\caption{Error estimator components $\errorest{\interp}{(t)},\errorest{\correction}{(t)},\errorest{\timestepping}{(t)}$ for the \rbl{parametric} double glazing problem { with $\nY=4$}}
	\label{fig:pde-l4-error-est-components-nointerp}
\end{figure}
The corresponding error estimator components are plotted in \cref{fig:pde-l4-error-est-components-nointerp}.
It can be seen that parametric refinement occurs when the interpolation error estimate $\errorest{\miset}$ reaches the dynamic tolerance $\interpTol$. In between these times, where the polynomial approximation space is \emph{fixed},
the interpolation error grows with time. The safety factor $\scalingInterpTol=10$ is visible in the value of $\interpTol$ 
relative to $\errorest{\miset,\letol}$. A consequence of  \rbl{including} $\scalingInterpTol$ is that the 
interpolation error is the dominant component \rbl{for the majority of the time integration}, 
therefore when we consider the error estimator \rbl{effectivity, what we are actually assessing is} the 
\emph{interpolation} error estimator. In fact, referring back to 
\cref{fig:pde-l4-error-eff-nointerp}, the error estimator efficiency is  \rbl{close to {one} apart from}
when the correction term $\errorest{\miset,\letol}$ \rbl{makes} a significant contribution (generally just after parametric refinement). 
\rbl{The resulting} overestimation is likely a consequence of the application of the triangle inequality in 
\rbl{\eqref{eq:correction-bound}}. 
The final observation is that {when} parametric refinement { is performed} we get a reduction of the interpolation error estimate, and this reduction is by approximately one order of magnitude. This corresponds to the choice of marking factor $\dorfler$.

\begin{figure}[!t]
	\begin{subfigure}{0.49\textwidth}
		\centering
		\includegraphics{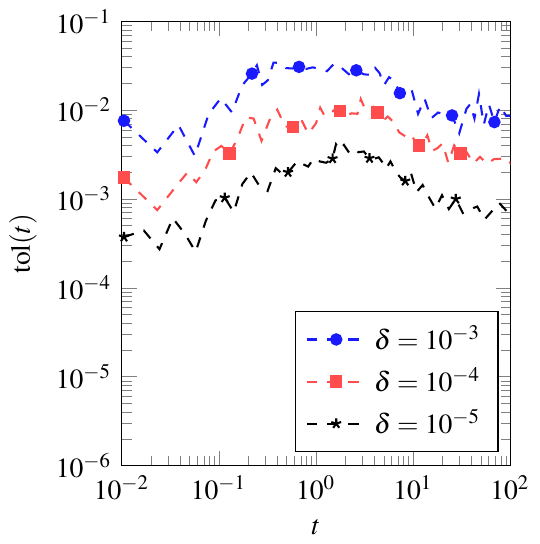}
		\caption{interpolation error tolerance $\interpTol(t)$}
		\label{fig:pde-l4-varied-letol-E}
	\end{subfigure}
	\begin{subfigure}{0.49\textwidth}
		\centering
		\includegraphics{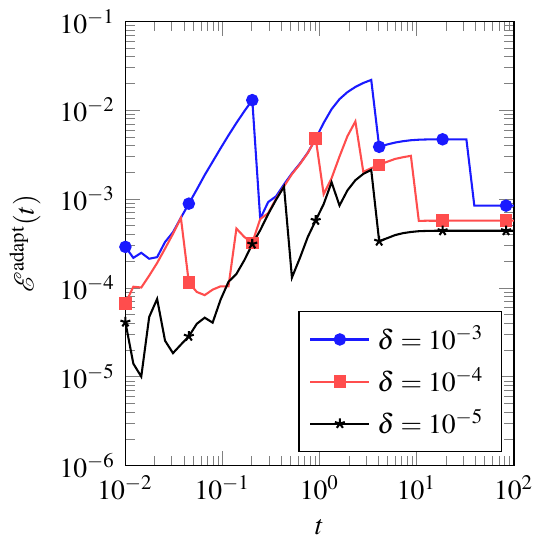}
		\caption{approximation error $\error{\adapt}(t)$}
		\label{fig:pde-l4-varied-letol-error}
	\end{subfigure}
	\caption{{\rbl{Interpolation error tolerance and associated approximation error evolution
			when varying the local error tolerance $\letol$} for the \rbl{parametric} double glazing problem with $\nY=4$}}  
	\label{fig:pde-l4-varied-letol}
\end{figure}

\rbl{We will now focus on the optimal selection of  the  inputs for \cref{alg:adaptive-loop} .}
\rbl{Reducing the local error tolerance $\letol$ will reduce} the global timestepping error for each collocation point. 
\rbl{The global timestepping error and the  local error tolerance are related through \eqref{eq:big-o-ge} with $\odeP=2$ when using \stabtr{}.}
\rbl{The definition of the correction estimator \eqref{eq:correction-est} suggests the interpolation error tolerance 
defined in \eqref{eq:toldef} will scale in the same way, so that $\interpTol \propto \letol^{2/3}$.} 
This \rbl{relationship} is demonstrated in \cref{fig:pde-l4-varied-letol-E}. 
\rbl{The plot in \cref{fig:pde-l4-varied-letol-error}  shows that the resulting approximation error} does not scale in the same way for all times. 
We note that the upper bound $\interpTol$ is briefly attained prior to the approximation being refined.

\begin{figure}[!ht]
	\begin{subfigure}{0.49\textwidth}
		\centering
		\includegraphics{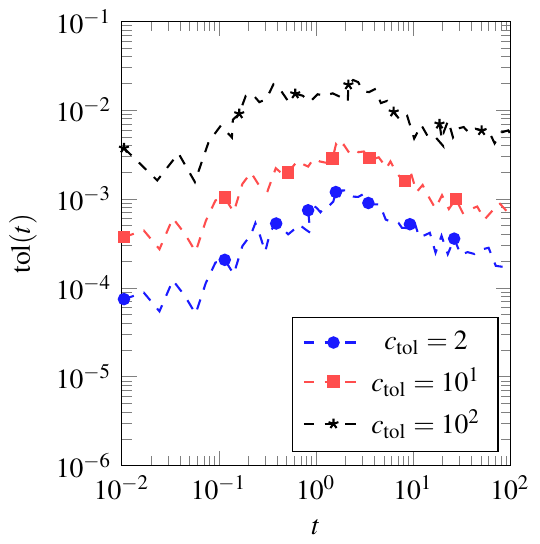}
		\caption{interpolation error tolerance $\interpTol(t)$}
		\label{fig:pde-l4-varied-safety-E}
	\end{subfigure}
	\begin{subfigure}{0.49\textwidth}
		\centering
		\includegraphics{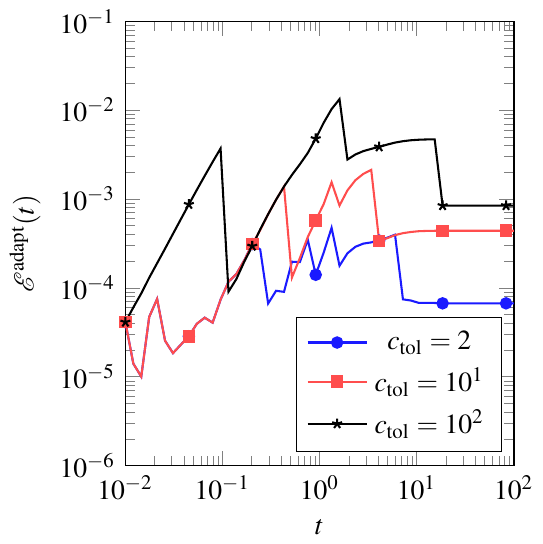}
		\caption{approximation error $\error{\adapt}(t)$}
		\label{fig:pde-l4-varied-safety-error}
	\end{subfigure}
	\caption{\rbl{Interpolation error tolerance  and associated approximation error evolution
			when varying the safety factor $\scalingInterpTol$} {for the \rbl{parametric} double glazing problem with $\nY=4$}}
	\label{fig:pde-l4-varied-safety}
\end{figure}

Modifying the safety factor $\scalingInterpTol$ affects the approximation by changing the relationship between $\interpTol$ 
and the correction estimator $\pi_{\miset,\letol}$. A smaller choice of $\scalingInterpTol$ should give an approximation 
with a smaller interpolation error. \Cref{fig:pde-l4-varied-safety-E} demonstrates the change in \rbl{the interpolation error tolerance}  as $\scalingInterpTol$ is varied. The resulting approximation errors 
are plotted in~\cref{fig:pde-l4-varied-safety-error}. 
We observe that the distinct peaks in the approximation error evolution are suppressed when the tolerance is tightened.
Note that a minimum value of $\scalingInterpTol=2$ is plotted because with $\scalingInterpTol=1$ the interpolation error estimates begin to be 
limited by timestepping errors. A sensible choice is   $\scalingInterpTol=10$.

\begin{figure}[!ht]			
	\includegraphics{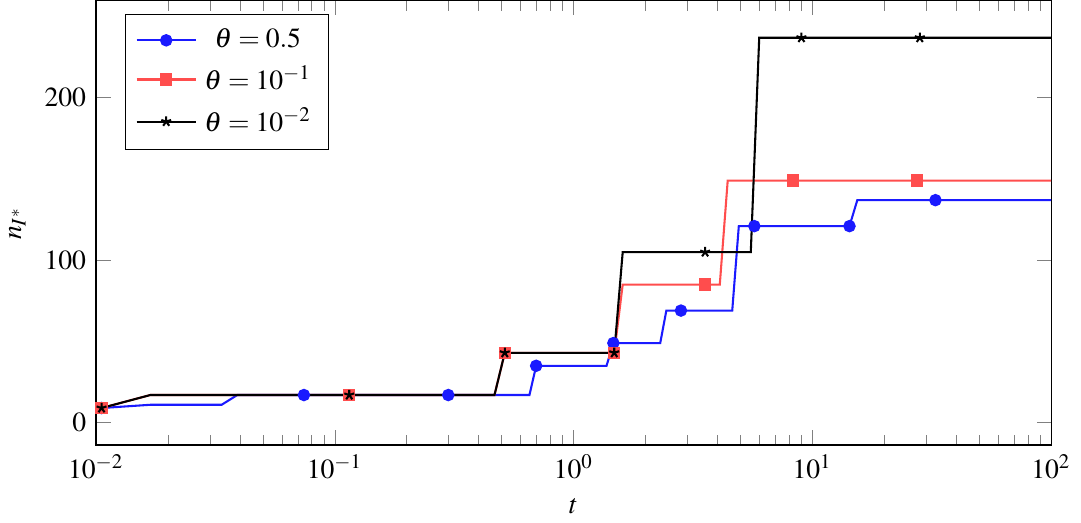}
	\caption{Number of collocation points {$\nColloc{\enhanced{\miset}}$} \rbl{needed for the approximation and error estimation} 
		as a function of time \rbl{when varying the marking parameter $\dorfler$} {for the \rbl{parametric} double glazing problem with $\nY=4$}} 
	\label{fig:ncolloc-marking}
\end{figure}

\rbl{Reducing the marking parameter $\dorfler$ 
	will  reduce the approximation error by a larger amount so the refinements will be less
	frequent. A large marking factor will give  a more prudent selection of marked multi-indices.}
\rbl{The results in} \Cref{fig:ncolloc-marking} demonstrate
that a more carefully constructed  multi-index set, corresponding to a larger $\dorfler$, generally results in 
fewer collocation points  being required which is computationally advantageous. 
The disadvantage is that more refinement steps \rbl{are needed}.  \rbl{Choosing $\dorfler=0.1$ is a reasonable compromise.}

{
\subsection{Exploiting stability for computational savings} \label{sub:interp}
We now}  discuss how the smoothing in time that is inherent in the solution of parabolic PDE problems can be exploited  to give computational savings.
By way of  motivation, consider two ODE systems of the form \eqref{eq:system-of-odes} for a particular parameter $\vCollocPt\in\Gamma$ with {initial} data $\vec{u}^{\h}(t_\iTime,\vCollocPt)$ and $\vec{v}^{\h}(t_\iTime,\vCollocPt)$ respectively at a time $t_\iTime \in [0,T]$.
The \rbl{respective solutions} $\vec{u}^{\h}(t,\vCollocPt)$ and $\vec{v}^{\h}(t,\vCollocPt)$ {for $t\in[t_{\iTime},T]$} can be represented as matrix exponentials and the difference between the two solutions $\vec{r}^{\h}(t,\vCollocPt)=\vec{u}^{\h}(t,\vCollocPt) - \vec{v}^{\h}(t,\vCollocPt)$ is given by
\begin{equation}
	\begin{aligned}
		\vec{r}^{\h}(t,\vCollocPt)
		=& 
		\exp\left(- (\massMatrix)^{-1} (\diffcnst \diffMatrix + \advMatrix(\vCollocPt) ) (t-t_\iTime)\right) \vec{r}^{\h}(t_\iTime,\vCollocPt).
	\end{aligned}
\end{equation}
If $\vec{r}^{\h}(t_\iTime,\vCollocPt)$ represents {the} error at time $t_\iTime$, then the behaviour of the norm of this error is controlled by the real components of the eigenvalues of the matrix $(\massMatrix)^{-1}\left(\diffcnst \diffMatrix + \advMatrix(\vCollocPt)\right)$. 
The diffusion matrix $\diffMatrix$ is symmetric and the advection matrix $\advMatrix(\vCollocPt)$ is skew-symmetric for all $\vCollocPt\in\Gamma$ hence the range of the real components of the eigenvalues {depends only on }the diffusion matrix~\cite{Loghin2006}. 
The diffusion matrix $\diffMatrix$ is positive-definite hence the eigenvalues are real and positive. 
\rbl{This means that,} for some constants $C_1,C_2 >0$ and $t\geq t_\iTime$
\begin{equation}
	\Vert \vec{r}^{\h}(t,\vCollocPt) \Vert_{\normX} \leq C_1 \exp(-C_2 \diffcnst(t-t_\iTime)).
	\label{eq:error-decay}
\end{equation}
We can interpret \eqref{eq:error-decay} as stating that an error at time $t_\iTime$ will (eventually) decay to zero through time. 
Note that there is no parameter dependence in the diffusion matrix, \rbl{so the bound is valid for all $\vCollocPt\in\Gamma$}.

\rbl{When solving the ODE system {\eqref{eq:system-of-odes}} using adaptive timestepping,} the accumulation of local error through time is related to the decay \eqref{eq:error-decay}. 
The local error made at each step will eventually decay through time and this is why local error control can give good global error {control}.
The bound \eqref{eq:error-decay} is independent of $\vCollocPt$ which suggests that the accumulation of local error into global error may be similar for all collocation points. 
This is in fact the case, as demonstrated in \cref{fig:pde-global-error-est}. 
\rbl{This insight} suggests that it should be possible to use a single global timestepping error estimate $\errorest{\gerror}(t) \approx \error{\h,\letol}(t;\vCollocPt)$ 
for all collocation points $\vCollocPt \in \collocSet{\adapt}$ \rbl{as a mechanism for reducing the overall} computational cost.
We could choose a fixed $\enhanced{\vCollocPt} \in \Gamma$, assume $\error{\h,\letol}(t;\enhanced{\vCollocPt}) \approx \error{\h,\letol}(t;\vCollocPt)$ for all $\vCollocPt \in \collocSet{\enhanced{\miset}}$ and then let $\errorest{\gerror}(t):=\errorest{\gerror}(t;\enhanced{\vCollocPt})$.
\rbl{{Since} the parametrisation of the  wind  in 
{\eqref{eq:eddies-wind}} gives $\Exp{\wind(\x,\y)} = \wind(\x,\vec{0})$,
in this case a sensible choice would be}
\begin{equation}
	\errorest{\gerror}(t) := \errorest{\gerror}(t;\vec{0}).
	\label{eq:error-est-mean-pt}
\end{equation}

\rbl{Whenever a new collocation point $\vCollocPt$ is added at  time $t$ in \cref{alg:adaptive-loop},}
we need access to the approximation $u^{\h,\letol}(\x,t;\vCollocPt)$. 
In \cref{ln:refine-initialise,ln:refine-catchup,ln:refine-update} this is achieved by applying the 
timestepping algorithm on the interval $[0,t]$ with the appropriate initial condition at $t=0$.
Alternatively, we could simply evaluate the current approximation $u^{\adapt}(\x,t,\vCollocPt)$ at time $t$ 
\rbl{at} the new collocation point $\vCollocPt$.
\rbl{To implement this strategy, we  simply need to} replace \cref{ln:refine-initialise,ln:refine-catchup,ln:refine-update} by
\begin{equation}
	u^{\h,\letol}_{local}(\x,t;\vCollocPt) \gets u^{\adapt}(\x,t,\vCollocPt) \label{eq:adapt-interp}
\end{equation}
with the corresponding timestep initialised by $\timestep_\vCollocPt \gets \timestep_0$. 
This \rbl{modification} 
 is again justified by the decay property \eqref{eq:error-decay}: 
the interpolation error introduced will decay through time.  \rbl{A significant computational saving will clearly result if this is done, because we do not have to solve {ODE systems} from time $t=0$ for {newly introduced} collocation points.}

\begin{figure}[!ht]
	\begin{subfigure}[t]{0.49\textwidth}
		\centering
		\includegraphics{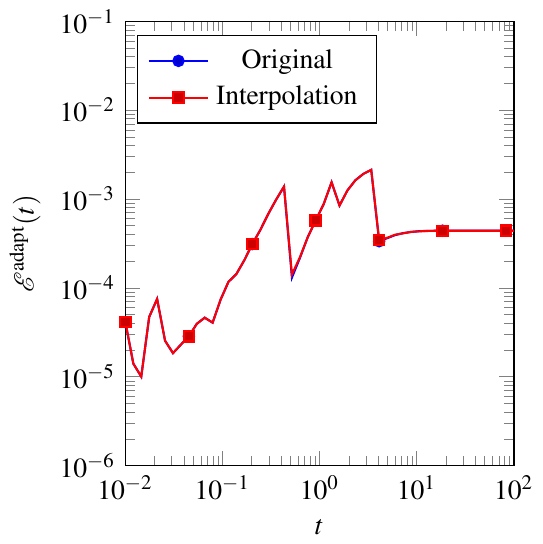}
		\caption{approximation error}
		\label{fig:initialise-interp-error}
	\end{subfigure}
	\begin{subfigure}[t]{0.49\textwidth}
		\centering
		\includegraphics{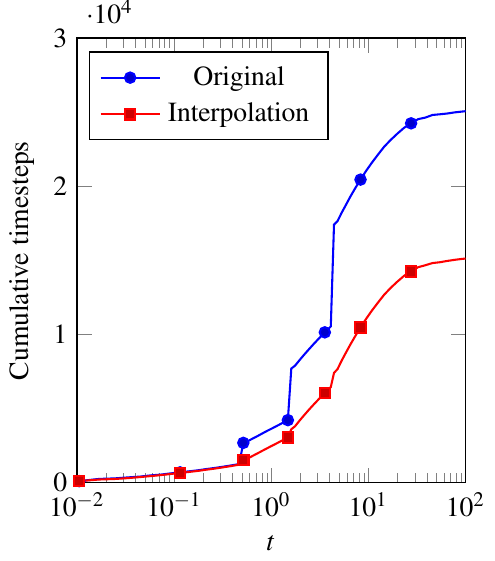}
		\caption{total timesteps \rbl{(all collocation points)}}
		\label{fig:initialise-interp-nsteps}
	\end{subfigure}
	\caption{Comparison of the performance of the original version of \cref{alg:adaptive-loop} and the modified version obtained by replacing lines 21-23 with the interpolation step \eqref{eq:adapt-interp}, {for} the parametric double glazing problem {with $\nY=4$}}
	\label{fig:initialise-interp}
\end{figure}

\rbl{The results in \cref{fig:initialise-interp-error} for this modified strategy demonstrate that initialisation through interpolation gives an approximation error profile that is essentially identical to that generated using \cref{alg:adaptive-loop}.}
If the computational cost is measured by counting the number of timesteps {used to approximate the ODE systems associated with} all the  collocation points, \rbl{then the results in \cref{fig:initialise-interp-nsteps} confirm that an equally accurate approximation has been constructed with a  significantly reduced  cost.}

{

\section{Computational Experiment: $\nY=64$}
\label{sub:example-sixtyfour}
We now consider a problem with a higher parametric dimension and investigate the adaptively constructed multi-index set $\miset(t)$ generated by \cref{alg:adaptive-loop}.
The test problem is set up as in \cref{sec:numerical-results} but now the wind field is constructed using a parametric stream function.
The method of construction is similar to that used in \cite{Phillips2015,Elman2018} and ensures $\wind(\x,\y)$ is divergence-free.
Specifically, $\wind(\x,\y) = \nabla \times \streamfn(\x,\y)$ where
\begin{equation}
	\streamfn(\x,\y) = \streamfn_0(\x) + \sum_{\iY=1}^{\nY} \sqrt{\lambda_\iY} \streamfn_{\iY}(\x) y_{\iY} \label{eq:stream-fn}
\end{equation}
and we choose $\y\in[-1,1]^{\nY}$ and $\rho(\y)=0.5^{\nY}$. Here $\{(\lambda_{\iY}, \streamfn_{\iY}(\x))\}_{\iY=1}^{\nY}$ are chosen to be the first $\nY$ eigenpairs (ordered in terms of decreasing $\lambda_{\iY})$ of an integral operator $\mathcal{C}:L^2(\spatialdomain)\to L^2(\spatialdomain)$ associated with a covariance function $C:\spatialdomain \times \spatialdomain \to\R$. That is, 
\begin{equation}
	\left(\mathcal{C}\streamfn_{\iY}\right)(\x):= \int_{\spatialdomain} C(\x,\x') \streamfn_{\iY}(\x')\dd \x' = \lambda_{\iY} \streamfn_{\iY}(\x). \label{eq:operatorC}
\end{equation}
Note that \eqref{eq:stream-fn} takes the form of a truncated Karhunen--Lo\'{e}ve expansion \cite[Section 5.4]{Lord2014}.
We choose the mean $\streamfn_{0}(\x):=-(1-x_1^2)(1-x_2^2)$ so that $\wind_0(\x) = \nabla \times \streamfn_0(\x)$ is equal to the mean field defined in \eqref{eq:w0_def} and we choose the covariance function
\begin{equation}
	C(\x_1,\x_2):= \left(\prod_{i,j=1}^{2}(1-x_{i,j}^2) \right) \sigma_0^2 \exp\left(\frac{\Vert\x_1 - \x_2\Vert_2^2}{L^2}\right) \label{eq:C}
\end{equation}
with $\sigma_0^2=5$ and $L=1$.  This is a modified version of the so-called Gaussian covariance function. The first term ensures that the stream function has zero variance on the boundary and therefore $\streamfn(\x,\y)=0$ for all $(\x,\y)\in\spatialbdry\times\Gamma$. The ensures that the wind field is always parallel to the boundary of $D$ and there is no inflow or outflow.

{
The required eigenpairs are approximated using collocation and quadrature on a suitably fine spatial grid \cite[Section 7.4]{Lord2014}.
The 2-norm of the wind perturbations \begin{equation}
	\wind_{\iY}(\x):=\sqrt{\lambda_\iY}\left(\nabla \times \streamfn_{\iY}(\x)\right)\label{eq:wi-rf}
\end{equation}
is computed on the same fine spatial grid and the maximum value over the grid points is plotted for each $\iY=1,2,...,64$ in \cref{fig:rf-eigs}.
The higher parameter dimensions should have a negligible effect on the solution as they correspond to very small perturbations. 
In this test problem we retain $\nY=64$ terms in the expansion of the stream function and allow \cref{alg:adaptive-loop} to identify a suitable multi-index set. 

}
We then use \cref{alg:adaptive-loop} to approximate the solution of the corresponding parametric advection--diffusion problem.
The input options used for the experiment are summarised in \cref{tab:ifissAndStabtr64} and we select the safety factor $\scalingInterpTol=10$ and marking parameter $\dorfler=10^{-1}$.
Approximations corresponding to newly introduced collocation points are added through interpolation as specified in \eqref{eq:adapt-interp}.
\begin{figure}
	\includegraphics{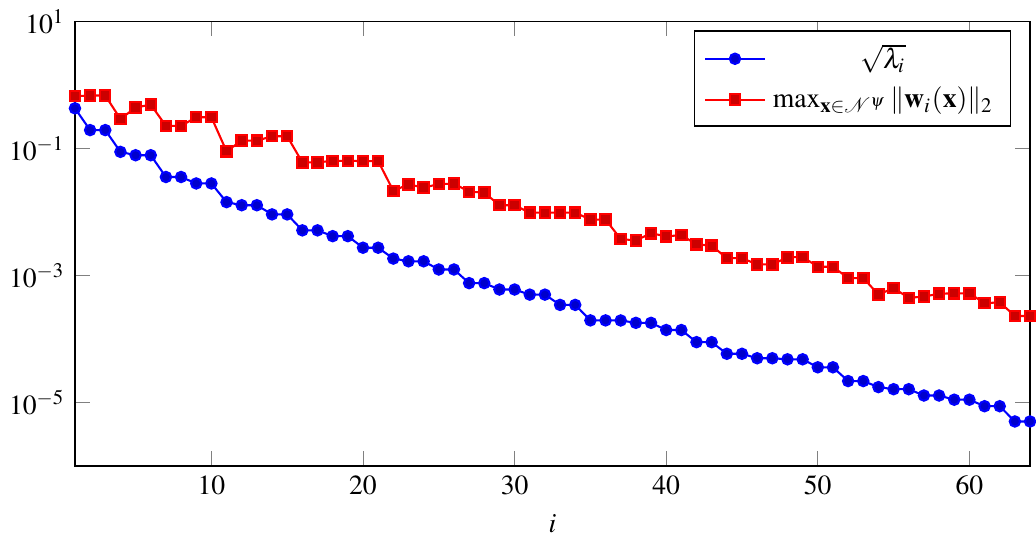}
	\caption{Approximations of the first $\nY=64$ eigenvalues defined by \eqref{eq:operatorC} and $\max_{\x\in\nodesRF} \Vert \wind_{\iY}(\x) \Vert_2$ where $\wind_{\iY}$ is as defined in \eqref{eq:wi-rf} and $\nodesRF$ are the spatial nodes used to approximate the eigenpairs}
	\label{fig:rf-eigs}
\end{figure}

The behaviour of the mean of the approximation is very similar to that obtained for the $d=4$ parameter problem (see \cref{fig:exp-pde-l4}) so is not shown.
This is not surprising as the same mean wind field was chosen in both cases. A more complicated space--time dependence is seen in the standard deviation in \cref{fig:stddev-pde-rf}. The error estimator and its components are plotted in \cref{fig:pde-l6-ests-rf}, demonstrating the control \cref{alg:adaptive-loop} provides for the interpolation error. The maximum level numbers selected in each dimension to construct the sparse grid interpolant, that is $\mi_{\iY}^{\max}:=\max_{\vMi\in\miset(t)} \mi_{\iY}$ for each $\iY=1,...,64$, at four stages of refinement, are shown in \cref{fig:maxmi-rf}.
We see that the maximum level number required grows through time and that \cref{alg:adaptive-loop} generally chooses higher levels in the parameter dimensions associated with larger eigenvalues $\lambda_{i}$ (corresponding to small $\iY$). Intuitively, we expect this as the lower dimensions correspond to stronger perturbations from the mean field.
{The final approximation is constructed using $1371$ collocation points. An additional $5880$ collocation points are used for the error estimator.
}

\begin{table}
	\caption{\ifiss{} and \stabtr{} input options used for approximating the parametric double glazing problem with $\nY=64$}
	\label{tab:ifissAndStabtr64}
	\begin{tabular}[t]{ll}
		\hline\noalign{\smallskip}
		\ifiss{} option & value\\
		\noalign{\smallskip}\hline\noalign{\smallskip}
		Hot wall time constant $\bdryrate$ & $10^{-1}$\\
		Grid parameter $\gridParam$  ($2^{2\gridParam}$ elements)& 6\\ 
		Stretched grid & no\\
		\rbl{Viscosity parameter $\epsilon$} & \rbl{$10^{-1}$ }\\
		\noalign{\smallskip}\hline
	\end{tabular} \qquad
	\begin{tabular}[t]{ll}
		\hline\noalign{\smallskip}
		\stabtr{} option & value\\
		\noalign{\smallskip}\hline\noalign{\smallskip}
		Averaging counter $\nstar$ & no averaging \\
		Initial timestep $\timestep_0$ & $10^{-9}$\\
		Local error tolerance $\letol$ & $10^{-6}$\\
		Final time $T$ & $100$ \\
		\noalign{\smallskip}\hline
	\end{tabular}
\end{table}

\begin{figure}[t]
	\begin{subfigure}{0.49\textwidth}
		\centering
		\includegraphics{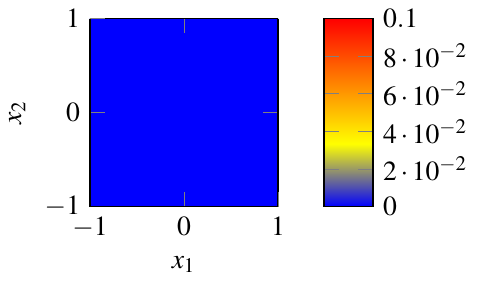}
		\caption{$\quad t\approx0.1$}
	\end{subfigure}\hspace{-10pt}
	\begin{subfigure}{0.49\textwidth}
		\centering
		\includegraphics{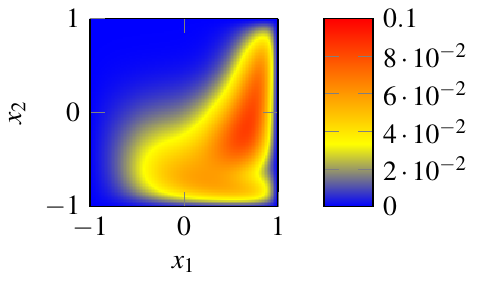}
		\caption{$\quad  t\approx1$}
	\end{subfigure}
	\begin{subfigure}{0.49\textwidth}
		\centering
		\includegraphics{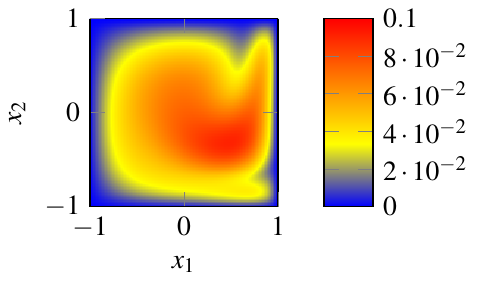}				
		\caption{$\quad t\approx10$}
	\end{subfigure}\hspace{-10pt}
	\begin{subfigure}{0.49\textwidth}
		\centering
		\includegraphics{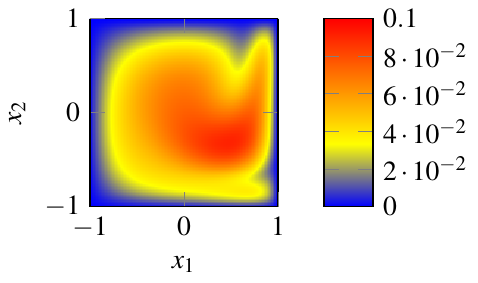}
		\caption{$\quad t\approx100$}
	\end{subfigure}
	\caption{\rbl{Snapshots in time of the standard deviation of the approximation to the solution of the parametric double glazing problem with $\nY=64$}}
	\label{fig:stddev-pde-rf}
\end{figure}


\begin{figure}
	\includegraphics{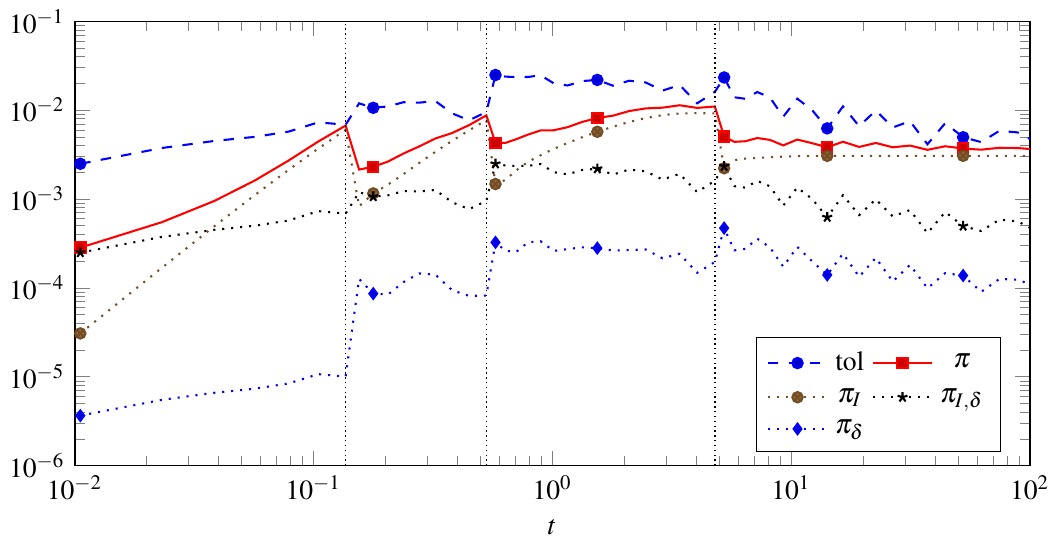}
	\caption{Error estimator $\errorest{}(t)$ and components $\errorest{\interp}(t),\errorest{\correction}(t),\errorest{\timestepping}(t)$ for the parametric double glazing problem with $\nY=64$. Vertical dotted lines denote \rbl{the parametric refinement} times.}
	\label{fig:pde-l6-ests-rf}
\end{figure}

\begin{figure}
	\includegraphics{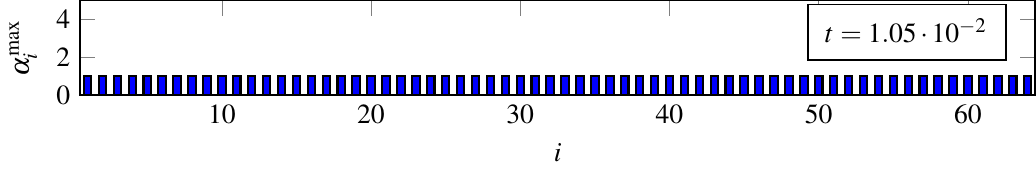}		\includegraphics{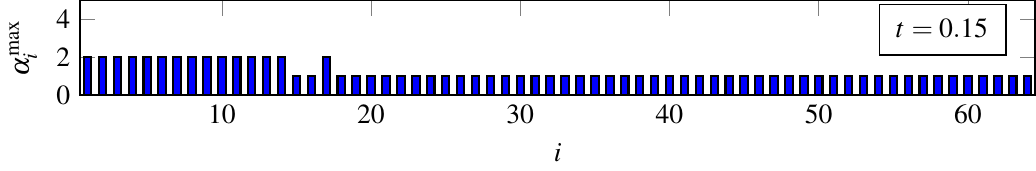}		\includegraphics{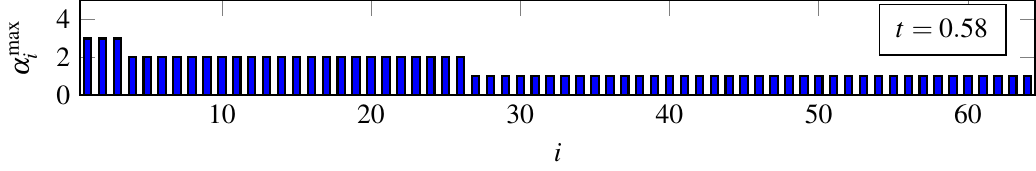}		\includegraphics{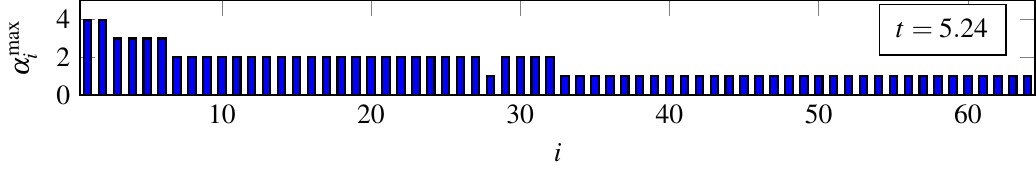}	
	\caption{Maximum level number $\mi_{\iY}^{\max}:=\max_{\vMi\in\miset(t)} \mi_{\iY}$ at four snapshots in time ($t\approx0.01,0.15,0.58,5.24$) for the parametric double glazing problem with $\nY=64$}
	\label{fig:maxmi-rf}
\end{figure}

Finally, in \cref{fig:nsteps-64} we illustrate the computational cost of applying \cref{alg:adaptive-loop}. Again this is measured as the total number of timesteps taken by the timestepping algorithm to solve all the ODE systems associated with all collocation points that are used in the computation of both the approximation and the error estimator. For short times, the bulk of the cost is in computing the interpolation error estimator. This is a fundamental issue with hierarchical error estimation and the doubling rule \eqref{eq:cc-doubling}. However, over time, as parametric enrichment is performed, numerical solutions to ODE systems associated with collocation points used at earlier times to compute the error estimator can be re-used to compute the approximation. In the long time limit we then see similar computational costs for the approximation and error estimator. 



 
To highlight the advantage of constructing a parametric approximation adaptively, we consider using a fixed naively chosen Smolyak sparse grid with $\nColloc{I_{SM}}=8321$ points corresponding to the multi-index set $\miset_{SM}:=\{\vMi \st \Vert \vMi \Vert_1 \leq 64 + 2\}$.
{This approximation would be exact for all polynomials of total degree less than or equal to 2 \cite[Theorem 4]{Barthelmann2000}.}
The cost of obtaining the associated approximation can be estimated by multiplying $\nColloc{I_{SM}}$ by the number of timesteps used by the adaptive timestepping method to solve the ODE system corresponding to the single collocation point $\y=\vec{0}$. As we can see in \cref{fig:nsteps-64}, constructing the approximation with the multi-index set $\miset_{SM}$ is more expensive than with our adaptive scheme. Moreover, since the multi-index set constructed by \cref{alg:adaptive-loop} contains higher level numbers, the associated approximation space contains polynomials of higher total degree than the one associated with $\miset_{SM}$.  A naive tensor grid approximation  in $d=4$ dimensions would be impossible.

\begin{figure}
	\includegraphics{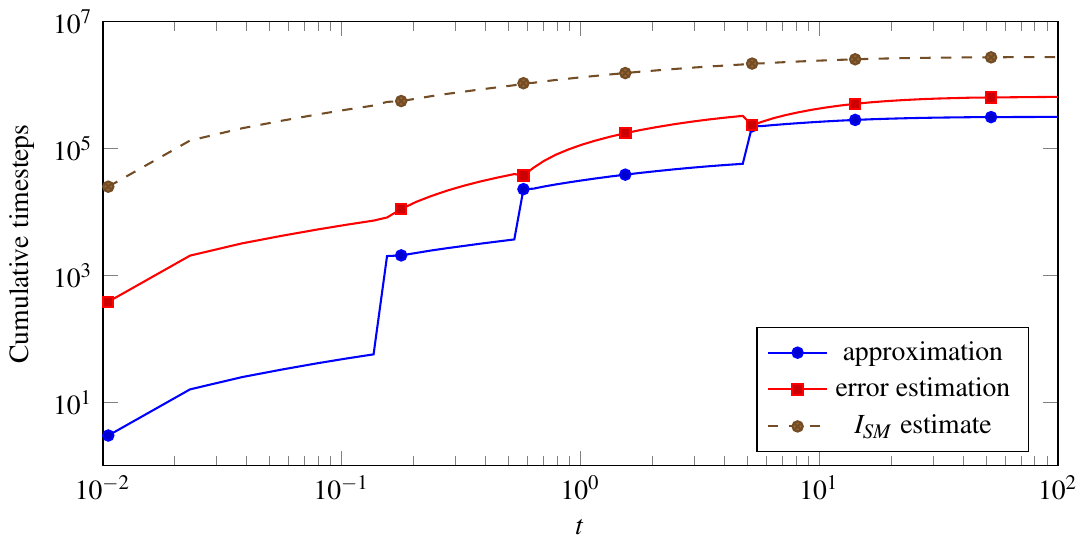}
	\caption{Total number of timesteps used in constructing the approximation, error estimator and a naive sparse grid approximation with multi-index set $\miset_{SM} = \{\vMi \st \Vert\vMi\Vert_{1} \leq 64 + 2\}$ for the parametric double glazing problem with $\nY=64$}
	\label{fig:nsteps-64}
\end{figure}


\section{Conclusions}
We introduced an adaptive nonintrusive stochastic collocation approximation scheme for a 
time-dependent, parametric PDE and  \rbl{developed} a corresponding error estimation strategy. 
The algorithm is used to monitor errors relating to the timestepping and construct an adaptive approximation
in which the interpolation error is controlled to a related tolerance.  \rbl{The efficiency of the error estimation strategy is demonstrated numerically, highlighting the necessity of adaptive stochastic collocation approximation}.  The adaptive approximation is able to control the error by refining the polynomial approximation space through time and sequentially constructing a \rbl{tailored approximation capturing the evolving uncertainty in the solution}.

The results presented demonstrate the control the algorithm inputs give over constructing an approximation. 
It is clear that tighter error control can be achieved by reducing the timestepping local error tolerance, but global error {control} in the {timestepping} would be required if we want rigorous control of the overall approximation error.
Stability of the underlying PDE plays an important role, firstly by ensuring local error control in the timestepping results in a good approximation and secondly by allowing the computational cost of constructing an approximation to be reduced through interpolation. 
The results can be shown to extend to higher dimensional parameter spaces.

\FloatBarrier
\bibliographystyle{spmpsci}      
\bibliography{bib}   



\section*{Statements \& Declarations}
\subsection*{Funding}
\thanks{{This work was partially supported by the EPSRC grants EP/V048376/1 and EP/W033801/1.}
%
\subsection*{Conflict of interest}
The authors declare that they have no conflict of interest.

\subsection*{Data Availability}
All data used in the manuscript is numerically generated using open source software packages \sparsekit{} and \ifiss{}. The MATLAB source code used to generate the numerical results will be made available on GitHub upon acceptance of the manuscript.

\end{document}